\documentclass[12pt]{amsart}

\usepackage{amsmath,amsfonts,amssymb,amscd,verbatim,delarray,verbatim}
  \textwidth
12.5cm \oddsidemargin 1cm \evensidemargin 1cm \textheight 19.5cm

\newcommand{\F}{{\mathbb F}}

\newcommand{\Q}{{\mathbb Q}}

\newcommand{\ch}{\operatorname{char}}

\newcommand{\di}{\operatorname{diag}}
\newcommand{\Imm}{\operatorname{Im}}

\DeclareMathOperator{\Mod}{mod}

\begin{document}
\numberwithin{equation}{section}

\newtheorem{theorem}{Theorem}[section]
\newtheorem{lemma}[theorem]{Lemma}

\newtheorem{prop}[theorem]{Proposition}
\newtheorem{proposition}[theorem]{Proposition}
\newtheorem{corollary}[theorem]{Corollary}
\newtheorem{corol}[theorem]{Corollary}
\newtheorem{conj}[theorem]{Conjecture}
\newtheorem{sublemma}[theorem]{Sublemma}

\theoremstyle{definition}
\newtheorem{defn}[theorem]{Definition}
\newtheorem{example}[theorem]{Example}
\newtheorem{examples}[theorem]{Examples}
\newtheorem{remarks}[theorem]{Remarks}
\newtheorem{remark}[theorem]{Remark}
\newtheorem{algorithm}[theorem]{Algorithm}
\newtheorem{question}[theorem]{Question}
\newtheorem{problem}[theorem]{Problem}
\newtheorem{subsec}[theorem]{}
\newtheorem{clai}[theorem]{Claim}

\def\toeq{{\stackrel{\sim}{\longrightarrow}}}
\def\into{{\hookrightarrow}}


\def\alp{{\alpha}}  \def\bet{{\beta}} \def\gam{{\gamma}}
 \def\del{{\delta}}
\def\eps{{\varepsilon}}
\def\kap{{\kappa}}                   \def\Chi{\text{X}}
\def\lam{{\lambda}}
 \def\sig{{\sigma}}  \def\vphi{{\varphi}} \def\om{{\omega}}
\def\Gam{{\Gamma}}   \def\Del{{\Delta}}
\def\Sig{{\Sigma}}   \def\Om{{\Omega}}
\def\ups{{\upsilon}}


\def\F{{\mathbb{F}}}
\def\BF{{\mathbb{F}}}
\def\BN{{\mathbb{N}}}
\def\Q{{\mathbb{Q}}}
\def\Ql{{\overline{\Q }_{\ell }}}
\def\CC{{\mathbb{C}}}
\def\R{{\mathbb R}}
\def\V{{\mathbf V}}
\def\D{{\mathbf D}}
\def\BZ{{\mathbb Z}}
\def\K{{\mathbf K}}
\def\XX{\mathbf{X}^*}
\def\xx{\mathbf{X}_*}

\def\AA{\Bbb A}
\def\BA{\mathbb A}
\def\HH{\mathbb H}
\def\PP{\Bbb P}

\def\Gm{{{\mathbb G}_{\textrm{m}}}}
\def\Gmk{{{\mathbb G}_{\textrm m,k}}}
\def\GmL{{\mathbb G_{{\textrm m},L}}}
\def\Ga{{{\mathbb G}_a}}

\def\Fb{{\overline{\F }}}
\def\Kb{{\overline K}}
\def\Yb{{\overline Y}}
\def\Xb{{\overline X}}
\def\Tb{{\overline T}}
\def\Bb{{\overline B}}
\def\Gb{{\bar{G}}}
\def\Ub{{\overline U}}
\def\Vb{{\overline V}}
\def\Hb{{\bar{H}}}
\def\kb{{\bar{k}}}

\def\Th{{\hat T}}
\def\Bh{{\hat B}}
\def\Gh{{\hat G}}

\def\cF{{\mathfrak{F}}}
\def\cC{{\mathcal C}}
\def\cU{{\mathcal U}}

\def\Xt{{\widetilde X}}
\def\Gt{{\widetilde G}}

\def\gg{{\mathfrak g}}
\def\hh{{\mathfrak h}}
\def\lie{\mathfrak a}

\def\minus{^{-1}}

\def\GL{\textrm{GL}}            \def\Stab{\textrm{Stab}}
\def\Gal{\textrm{Gal}}          \def\Aut{\textrm{Aut\,}}
\def\Lie{\textrm{Lie\,}}        \def\Ext{\textrm{Ext}}
\def\PSL{\textrm{PSL}}          \def\SL{\textrm{SL}}
\def\loc{\textrm{loc}}
\def\coker{\textrm{coker\,}}    \def\Hom{\textrm{Hom}}
\def\im{\textrm{im\,}}           \def\int{\textrm{int}}
\def\inv{\textrm{inv}}           \def\can{\textrm{can}}
\def\id{\textrm{id}}              \def\Char{\textrm{char}}
\def\Cl{\textrm{Cl}}
\def\Sz{\textrm{Sz}}
\def\ad{\textrm{ad\,}}
\def\SU{\textrm{SU}}
\def\Sp{\textrm{Sp}}
\def\PSL{\textrm{PSL}}
\def\PSU{\textrm{PSU}}
\def\rk{\textrm{rk}}
\def\PGL{\textrm{PGL}}
\def\Ker{\textrm{Ker}}
\def\Ob{\textrm{Ob}}
\def\Var{\textrm{Var}}
\def\poSet{\textrm{poSet}}
\def\Al{\textrm{Al}}
\def\Int{\textrm{Int}}
\def\Mod{\textrm{Mod}}
\def\Smg{\textrm{Smg}}
\def\ISmg{\textrm{ISmg}}
\def\Ass{\textrm{Ass}}
\def\Grp{\textrm{Grp}}
\def\Com{\textrm{Com}}
\def\rank{\textrm{rank}}

\def\char{\textrm{char}}

\newcommand{\Or}{\operatorname{O}}

\def\tors{_\def{\textrm{tors}}}      \def\tor{^{\textrm{tor}}}
\def\red{^{\textrm{red}}}         \def\nt{^{\textrm{ssu}}}

\def\sss{^{\textrm{ss}}}          \def\uu{^{\textrm{u}}}
\def\mm{^{\textrm{m}}}
\def\tm{^\times}                  \def\mult{^{\textrm{mult}}}

\def\uss{^{\textrm{ssu}}}         \def\ssu{^{\textrm{ssu}}}
\def\comp{_{\textrm{c}}}
\def\ab{_{\textrm{ab}}}

\def\et{_{\textrm{\'et}}}
\def\nr{_{\textrm{nr}}}

\def\nil{_{\textrm{nil}}}
\def\sol{_{\textrm{sol}}}
\def\End{\textrm{End\,}}

\def\til{\;\widetilde{}\;}



\title[A commutator description of the solvable radical]
{{\bf A commutator description of the solvable radical of a finite
group  }}

\author[Gordeev, Grunewald, Kunyavskii,  Plotkin] {Nikolai
Gordeev, Fritz Grunewald,  Boris Kunyavskii, Eugene Plotkin }
\address{Gordeev: Department of Mathematics, Herzen State
Pedagogical University,
48 Moika Embankment, 191186, St.Petersburg, RUSSIA} \email{nickgordeev@mail.ru}

\address{Grunewald: Mathematisches Institut der
Heinrich-Heine-Universit\"at D\"usseldorf, Universit\"atsstr. 1, 40225
D\"usseldorf, GERMANY} \email{grunewald@math.uni-duesseldorf.de}

\address{ Kunyavskii: Department of
Mathematics, Bar-Ilan University, 52900 Ramat Gan, ISRAEL}
\email{kunyav@macs.biu.ac.il}
\address{ Plotkin: Department of
Mathematics, Bar-Ilan University, 52900 Ramat Gan, ISRAEL}
\email{plotkin@macs.biu.ac.il}


\begin{abstract}
We are looking for the smallest integer $k>1$ providing the
following characterization of the solvable radical $R(G)$ of any
finite group $G$: $R(G)$ coincides with the collection of $g\in G$
such that for any $k$ elements $a_1, a_2,\dots ,a_k \in G$ the
subgroup generated by the elements $g, a_iga_i^{-1}$, $i=1,\ldots
,k$, is solvable. We consider a similar problem of finding the
smallest integer $\ell >1$ with the property that $R(G)$ coincides
with the collection of $g\in G$ such that for any $\ell$ elements
$b_1, b_2,\dots ,b_\ell \in G$ the subgroup generated by the
commutators $[g,b_i]$, $i=1,\dots ,\ell$, is solvable.
Conjecturally, $k=\ell =3$. We prove that both $k$ and $\ell$ are
at most 7. In particular, this means that a finite group $G$ is
solvable if and only if every 8 conjugate elements of $G$ generate
a solvable subgroup.
\end{abstract}

\maketitle

\thispagestyle{empty} \vspace{1.0cm} 

\tableofcontents

\section{Introduction} \label{sec:intro}

\subsection{Main results}\label{subsec:mnres}

Let $F_2=F(x,y)$ be the free two generator group. Define a
sequence $ \overrightarrow{e}=e_1, e_2, e_3,\dots, $ where
$e_i(x,y)\in F_2$, by
$$
e_1(x,y)=[x,y]=xyx^{-1}y^{-1},\ldots,
e_n(x,y)=[e_{n-1}(x,y),y],\ldots,
$$

An element $g$ of a group $G$ is called an Engel element if for
every $a\in G$ there exists a number $n=n(a,g)$ such that
$e_n(a,g)=1$.

In 1957 R. Baer proved the following theorem \cite{Ba}, \cite{H}:

\begin{theorem}\label{th:br}
The nilpotent radical of a noetherian group $G$ coincides with the
collection of all Engel elements of $G$.
\end{theorem}

In particular, Baer's theorem is true for finite groups. Similar
theorems have been established for many classes of infinite groups
satisfying some additional conditions (see for example \cite{Plo},
\cite{Pla}).

A tempting but difficult problem  is to find a counterpart of
Baer's theorem for the solvable radical of a finite group, in
other words, to find an Engel-like sequence $\overrightarrow
u=u_n(x,y)$ such that an element $g$ of a finite group $G$ belongs
to the solvable radical $R(G)$ of $G$ if and only if for any $a\in
G$ there exists a number $n=n(a,g)$ such that $u_n(a,g)=1$. The
first results towards a solution of this problem have been
obtained in \cite{BGGKPP1}, \cite{BGGKPP2}, \cite{BWW}, and
\cite{BBGKP}.

In the paper \cite{GKPS} a Thompson-like characterization of the
solvable radical of  finite groups (and, more generally, linear
groups and PI-groups)  has been obtained.

\begin{theorem}\cite{GKPS}\label{rad:aner}
The solvable radical $R(G)$ of a finite group $G$ coincides with
the set of all elements $g\in G$ with the following property: for
any $a\in G$ the subgroup generated by $g$ and $a$ is solvable.
\end{theorem}

This theorem can be viewed as an implicit description of the
solvable radical since it does not provide any  explicit formulas
which determine if a particular element belongs to $R(G)$.

In the present paper our goal is to obtain a new characterization
of the solvable radical $R(G)$ of a finite group $G$.

\begin{theorem}\label{main:solv_2}
The solvable radical of any finite group $G$ coincides with the
collection of $g\in G$ satisfying the property: for any 7 elements
$a_1, a_2,\dots ,a_7 \in G$ the subgroup generated by the elements
$g, a_iga_i^{-1}$, $i=1,\dots ,7$, is solvable.
\end{theorem}

The proof involves the classification of finite simple groups.

This theorem implies the following characterization of finite
solvable groups:

\begin{theorem}\label{main:solv}
A finite group $G$ is solvable if and only if every $8$ conjugate
elements of $G$ generate a solvable subgroup.
\end{theorem}

We hope to sharpen these results.

\begin{conj} \label{conjecture:rad}
The solvable radical of a finite group $G$ coincides with the
collection of $g\in G$ satisfying the property: for any 3 elements
$a, b, c \in G$ the subgroup generated by the conjugates $g,
aga^{-1}, bgb^{-1}, cgc^{-1}$ is solvable.
\end{conj}

This statement implies

\begin{conj} \label{conjecture:solv}
A finite group $G$ is solvable if and only if every four conjugate
elements of $G$ generate a solvable subgroup.
\end{conj}

\begin{remark} \label{rem:sharp}
These characterizations are the best possible: in the symmetric
groups $S_n$ $(n\geq 5)$ any triple of transpositions generates a
solvable subgroup.
\end{remark}

\begin{remark} \label{rem:imply}
The main step in our proof of Theorem \ref{main:solv_2} is Theorem
\ref{th:radel_1} below. To prove Conjecture \ref{conjecture:rad}
(and hence Conjecture \ref{conjecture:solv}), one has to extend
the statement of Theorem \ref{th:radel_1} to all {\em almost
simple} groups, i.e. to the groups $H$ such that $G\subseteq
H\subseteq \Aut(G)$ for some simple group $G$.
\end{remark}

\begin{remark} \label{rem:linear}
The statements of Theorems \ref{main:solv_2} and \ref{main:solv}
remain true for arbitrary linear groups. Once Conjectures
\ref{conjecture:rad} and \ref{conjecture:solv} are proved, they
can also be extended to arbitrary linear groups.
\end{remark}


Throughout the paper $\langle a_1,\ldots, a_k\rangle$ stands for
the subgroup of $G$  generated by $a_1,\ldots, a_k\in G$. We
define the commutator of $x, y\in G$ by $[x,y]=xyx^{-1}y^{-1}$.

\begin{defn} \label{defn:radel} Let $k\ge2$ be an integer.
We say that $g\in G$ is a {\em $k$-radical element} if for any
$a_1,\dots ,a_k\in G$ the subgroup $H=\langle [a_1,g],\dots ,
[a_k,g]\rangle$ is solvable.
\end{defn}

We prove the following result.

\begin{theorem} \label{th:radel_1} Let
$G$ be a finite nonabelian simple group. Then $G$ does not contain
nontrivial $3$-radical elements.
\end{theorem}

This theorem implies  Theorems \ref{main:solv_2} and
\ref{main:solv}.

The proof goes by case-by-case inspection of simple groups
(alternating groups, groups of Lie type, sporadic groups). In fact
we  prove a more precise result (Theorem \ref{th:lie1}) which
distinguishes between 2-radical and 3-radical elements.

The following simple fact allows us to define a new invariant
of a finite group.

\begin{prop}
Let $G$ be a group which has no nontrivial solvable normal
subgroups. Then for every $g \in G, g \ne 1$ the group $H_g =
\left<[g, G]\right>$ is not solvable.
\end{prop}

\begin{proof}
For every $x, y \in G$ we have
$$[g,x]^{-1}[g,y] = (x g x^{-1}g^{-1})(g y g^{-1}y^{-1}) = (x g
x^{-1})(y g^{-1}y^{-1}) \in H_g.$$ Thus, $C_gC_{g^{-1}} \subset
H_g$ where $C_g, C_{g^{-1}}$ are the corresponding conjugacy
classes. Since the set $C_gC_{g^{-1}}$ is invariant under
conjugation, the subgroup $F = \left<C_gC_{g^{-1}}\right> \leq
H_g$ is normal in $G$ and therefore cannot be solvable.
\end{proof}

\begin{corollary} \label{cor:k-rad}
Let $G$ be a finite group, and let $R(G)$ denote the solvable
radical of $G$. Then $g \notin R(G)$ if and only if there exist an
integer $n$ and $x_1, \dots ,x_n\in G$ such that the subgroup
$\left<[g, x_1], \dots ,[g, x_n]\right>$ is not solvable.
\end{corollary}

\begin{defn}
Denote by $\kappa (g)$ the smallest possible $n$ with the
following property: $g \notin R(G)$ if and only if there exist
$x_1, \dots ,x_n\in G$ such that the subgroup $\left<[g, x_1],
\dots ,[g, x_n]\right>$ is not solvable. We call the number
$\kappa (G):=\max_{g\in G\setminus R(G)}\kappa (g)$ the radical
degree of $G$.
\end{defn}

In these terms we have to prove that the radical degree of a finite
nonabelian simple group $G$ is $\leq 3$. Our most precise result,
which implies Theorem \ref{th:radel_1} and, correspondingly,
Theorems \ref{main:solv_2} and \ref{main:solv}, is the following

\begin{theorem} \label{th:lie1}
If $G$ is a finite nonabelian simple group, then $\kappa (G)\leq
3$. If $G$ is a group of Lie type over a field $K$ with $\ch K
\neq 2$ and $K \neq {\mathbb F}_3$, or a sporadic group not
isomorphic to $Fi_{22}$ or $Fi_{23}$, then $\kappa (G)=2.$
\end{theorem}

\subsection{Notation and conventions}


First introduce some standard notation which mostly follows
\cite{St}, \cite{Ca1}, \cite{Ca2}.

Denote by $G=G(\Phi,K)$ a Chevalley group where $\Phi$ is a
reduced irreducible root system and $K$ is a field. Assume that
$\Phi$ is generated by a set of simple roots
$\Pi=\{\alpha_1,\dots,\alpha_r\}$, that is $\Phi=
\left<\alpha_1,\dots,\alpha_r\right>$. We number the roots
according to \cite{Bou}. Let $W=W(\Phi)$ be the Weyl group
corresponding to $\Phi$. Denote by $\Phi^+$, $\Phi^-$ the sets of
positive and negative roots, respectively.

We use the standard notation $u_\alpha(t)$, $\alpha\in \Phi$, $t\in
K$, for elementary root unipotent elements of $G$. Correspondingly,
split semisimple elements will be denoted by $h_\alpha(t),$ $t\in
K^*$, where $K^*$ is the multiplicative group of $K$. For $\alpha\in
\Phi$, let $U_\alpha$ denote the root subgroup generated by all
elementary root unipotent elements $u_\alpha(t)$.

For the sake of completeness, recall that
$w_\alpha(t)=u_\alpha(t)u_{-\alpha}(-t^{-1})u_\alpha(t)$,
$w_\alpha=w_\alpha(1)$ and $h_\alpha(t)=w_\alpha(t)w_\alpha^{-1}$.
Define the subgroups $U=U^+=\left<u_\alpha(t), \alpha\in \Phi^+,
t\in K\right>$, $V=U^-=\left<u_\alpha(t), \alpha\in \Phi^-, t\in
K\right>$, $T=\left<h_\alpha(t), \alpha\in \Phi, t\in K^*\right>$,
and $N=\left<w_\alpha(t), \alpha\in \Phi, t\in K^*\right>$.

As usual, the Borel subgroups $B^\pm$ are $B=B^+=TU$, $B^-=TU^-$.
The group $N$ contains $T$, and $N/T\cong W$.
Denote by $\dot w$ a preimage of $w\in W$ in $N$.

We also consider twisted Chevalley groups over finite fields.
Assume that $K$ is a finite field of characteristic $p$ and
$|K|=q=p^s$. By a twisted Chevalley group we mean the group
$G^F=G^F(\Phi,\overline K)$ of fixed points of the simply
connected Chevalley group $G(\Phi,\overline K)$ under the
Frobenius map $F$ (see \cite{St}, \cite{Ca1}, \cite{Ca2}). Here
$\overline K$ stands for the algebraic closure of $K$. Let
$\theta$ be the field automorphism corresponding to $F$. Denote by
$k=K^\theta$ the subfield of $\theta$-fixed points for all cases
except for the Suzuki groups and the Ree groups. For the latter
groups suppose that $k=K$. Let $\gamma$ be the graph automorphism
corresponding to $F$. We denote by $\Phi^\gamma$ the root system
which determines the structure of the group
$G^F=G^F(\Phi,\overline K)$. Elementary root unipotent elements
$u_\alpha(t)$, $\alpha \in \Phi^\gamma$, have either one parameter
$t\in K$ or $t\in k$, or two parameters $u_\alpha (t,u),$ $t,u\in
K$ (for the cases $^2A_2$, $^2B_2$, $^2F_4$), or three parameters
$u_\alpha (t,u,v),$ $t,u,v\in K$ (for $^2G_2$), see \cite{St}.
Again, the root subgroups $U_\alpha$ are generated by root
unipotent elements. The subgroups $B^F$, $W^F$, $T^F$, $H^F$,
${U^\pm}^F$ in $G^F$ are defined in a standard way, see
\cite{Ca2}. A maximal torus of $G^F$ is a subgroup of the form
$T^F$, where $T$ is an $F$-stable maximal torus of $G$.  A maximal
torus  $T^F$ is called quasisplit if it is contained in $B^F$.
Throughout the paper we suppress the map $F$ in the notations. We
also suppress $\gamma$ in the notation of the root system
corresponding to the group $G^F$. Whenever we need to specify the
type of a group, it will be written explicitly.

We follow \cite{Ca2} in the notation of twisted forms. Thus
unitary groups are denoted by $PSU_n(q^2)$ (and not by
$PSU_n(q)$), the notation $^2F_4(2^{2m+1})$ means that
$q=\sqrt{2^{2m+1}}$, etc.

The paper is organized as follows. In Section \ref{sec:reduction}
we reduce Theorem \ref{main:solv} to Theorem \ref{th:lie1}. In
Sections \ref{sec:alt}--\ref{sec:spor} we prove Theorem
\ref{th:lie1} using case-by-case analysis.

\noindent {\it Acknowledgements}. Gordeev was partially supported
by the INTAS grant N-05-1000008-8118. Kunyavski\u\i \ and Plotkin
were partially supported by the Ministry of Absorption (Israel),
the Israeli Science Foundation founded by the Israeli Academy of
Sciences --- Center of Excellence Program, the Minerva Foundation
through the Emmy Noether Research Institute of Mathematics, and by
the RTN network HPRN-CT-2002-00287. A substantial part of this
work was done during Gordeev's visits to Bar-Ilan University in
May 2005 and May 2006 (partially supported by the same RTN
network) and the visit of Grunewald, Kunyavski\u\i \ and Plotkin
to MPIM (Bonn) during the activity on ``Geometry and Group
Theory'' in July 2006.  The support of these institutions is
highly appreciated.

We are very grateful to J.~N.~Bray, B.~I.~Plotkin, and
N.~A.~Vavilov for useful discussions and correspondence. Our
thanks also go to the anonymous referee for numerous remarks.


\section{Reduction  Theorem  } \label{sec:reduction}

Let us show how Theorem \ref{th:lie1} implies Theorem
\ref{main:solv}.

Suppose Theorem \ref{th:lie1} is proven, and let us show that the
solvable radical $R(G)$ of a finite group $G$ coincides with the
collection of $g\in G$ satisfying the property: for any $7$
elements $a_1, a_2,\dots ,a_7 \in G$ the subgroup generated by the
elements $g, a_iga_i^{-1}$, $i=1,\dots ,7$, is solvable.

For the sake of convenience, let us call the elements $g\in G$
satisfying the condition of the theorem, suitable.

Suppose  $g\in R(G)$. Since $R(G)$ is a normal subgroup,
$aga^{-1}$ belongs to $R(G)$ for any $a\in G$. Hence for any $k$
the subgroup $\langle a_1ga_1^{-1},\ldots ,a_kga_k^{-1}\rangle $,
where $a_1,\dots ,a_k\in G$, is solvable. Therefore, all the
elements of $R(G)$ are suitable.

Suppose now that $g\in G$ is a suitable element. We want to show
that $g$ belongs to $R(G)$.  It is enough to prove that there are
no non-trivial suitable elements in the semisimple group $G/R(G)$.
So one can assume that the group $G$ is semisimple in the sense
that $R(G)=1$.

As usual we consider a minimal counterexample $G$ to the statement
above.

Recall that any finite semisimple group $G$ contains a unique
maximal normal centreless completely reducible (CR) subgroup (by
definition, CR means a direct product of finite non-abelian simple
groups) called the CR-radical of $G$ (see \cite[3.3.16]{Ro}). We
call a product of the isomorphic factors in the decomposition of
the $CR$-radical {\it an isotypic component} of $G$. Denote the
$CR$-radical of $G$ by $V$. This is a characteristic subgroup of
$G$.

Let us show that $V$ has only one isotypic component.  Suppose
$V=N_1\times N_2$, where $N_1\cap N_2=1$. Consider $\bar G=G/N_1$
and denote $\bar R=R(G/N_1)$. Consider a suitable $g\in G$, $g\neq
1$ and denote by $\bar g$ (resp. ${\bar{\bar g}}$) the image of
$g$ in $\bar G$ (resp. $\bar G/\bar R$). Since $\bar G/\bar R$ is
semisimple and $\bar{\bar g}\in\bar G/\bar R$ is suitable, we have
${\bar{\bar g}}=1$ (because $G$ is a minimal counter-example) and
hence $\bar g\in R(G/N_1)$. Consider $V/N_1\simeq N_2$. Then
$V/N_1\subset G/N_1$ is semisimple and therefore $V/N_1\cap
R(G/N_1)=1$. Since $\bar g\in R(G/N_1)$, we have $[\bar g, \bar
v]=\bar 1$ for every $\bar v \in V/N_1$. Hence $[g,v]\in N_1$ for
every $v\in V$. Similarly, $[g,v]\in N_2$ for every $v\in V$.
Therefore $[g,v]=1$. Hence $g$ centralizes every $v\in V$. Since
the centralizer of $V$ in $G$ is trivial, we get $g=1$.
Contradiction.

Any $g\in G$ acts as an automorphism $\tilde g$ on $V=H_1\times
\cdots\times H_n$, where  all $H_i,$ $1\leq\ i\leq n$, are
isomorphic nonabelian simple groups.

Suppose that $g$ is a suitable element. Let us show that $\tilde
g$ cannot act on $V$ as a non-identity element of the symmetric
group $S_n$. Denote by $\sigma$ the element of $S_n$ corresponding
to $\tilde g$.

By definition, the subgroup $\Gamma=\langle g,
x_igx_i^{-1}\rangle$, $i=1,\dots ,7$, is solvable for any elements
$x_i\in G$. Evidently, the subgroup $\langle [g,x_1],
[g,x_2]\rangle$ lies in $\Gamma$.

Suppose $\sigma\neq 1$, and so $\sigma (k)\ne k$ for some $k\leq
n$. Take $\bar x_1$ and $\bar x_2$ of the form $\bar x_i =
(1,\ldots,x_i^{(k)},\ldots, 1)$, where $x_i^{(k)}\neq 1$ lies in
$H_k$ $(i=1,2)$. Then we may assume $(\bar x_i)
^\sigma=(x_i^{(k)},1,\ldots,1)$, and so $[g,\bar x_i]=(\bar x_i)
^\sigma \bar
x_i^{-1}=(x_i^{(k)},1,\ldots,{(x_i^{(k)}})^{-1},\ldots, 1).$


By a theorem of Steinberg, $H_k$ is generated by two elements, say
$a$ and $b$. On setting $x_1^{(k)}=a$, $x_2^{(k)}=b$, we conclude
that the group generated by $[g,\bar x_1]$ and $[g,\bar x_2]$
cannot be solvable because the first components of these elements,
$a$ and $b$, generate the simple group $H_k$. Contradiction with
solvability of $\Gamma$.

So we can assume that a suitable element $g\in G$ acts as an
automorphism of a simple group $H$. Then we  consider the
extension of the group $H$ with the automorphism $\tilde g$.
Denote this almost simple  group by $G_1$.  We shall use the
formula
$$
y[x,g]y^{-1}=[x,g][[g,x],y].
$$

Since $G_1$  has no centre, one can choose $x\in H$ such that
$[x,\tilde g]\neq 1.$ Evidently, $[x,\tilde g]$ belongs to the
simple group $H$. Then by Theorem \ref{th:lie1}, there exist $y_1,
y_2, y_3\in H$ such that the subgroup $\langle [[x,g],y_1],
[[x,g],y_2],[[x,g],y_3]\rangle $ is not solvable.  But
$$
\langle [[x,g],y_1], [[x,g],y_2],[[x,g],y_3]\rangle  \leq \langle
y_i[x,g]y_i^{-1}, [x,g]\ |\ i=1, 2, 3\rangle \leq $$ $$\leq
\langle g, x^{-1}gx, y_i^{-1}gy_i, y_i^{-1}x^{-1}gxy_i\ |\
i=1,2,3\rangle
$$
Since $g$ is suitable, the latter subgroup must be solvable.
Contradiction with the choice of $y_i$.


\section{Alternating groups} \label{sec:alt}

\begin{prop}
Let $G =A_n$, $n\geq 5$. Then $\kappa (G) = 2$.
\end{prop}

\begin{proof}
For $n=5,6,7$ the statement can be checked in a straightforward
manner, so assume $n\geq 8$. Let us proceed by induction. Let
$y\in G$, $y\neq 1$. First suppose that $y$ can be written in the
form
\begin{equation} \label{eq:ind}
y=\sigma\tau , \quad \sigma\in A_m, \quad \sigma\neq 1, \quad
5\leq m<n,
\end{equation}
where $\sigma$ and $\tau$ are disjoint (and thus commute). Then by
induction hypothesis there exist $\sigma_1,\sigma_2\in A_m$ such
that the subgroup generated by $[\sigma,\sigma_1]$ and
$[\sigma,\sigma_2]$ is not solvable. Take $x_i=\sigma_i\tau$,
$i=1,2$. Then $[y,x_i]=[\sigma,\sigma_i]$, and we are done.

Suppose $y$ cannot be represented in the form (\ref{eq:ind}). Then
we have one of the following cases: either $n$ is odd and
$y=(12\dots n)$, or $n$ is even and $y=(12\dots n-2)(n-1,n)$. In
any of these cases we take $x_1=(123)$ and $x_2=(345)$ and get
$\left<[x_1,y],[x_2,y]\right>\cong A_5$.
\end{proof}

\section{Groups of Lie type of small rank over fields of large characteristic}\label
{sec:lie_1}

\begin{prop} \label{prop:rank1}
Let $G$ be one of the groups $A_1(q)$ $(q \ne 2, 3)$, $^2A_2(q^2)$
$(q \ne 2)$, $^2B_2(2^{2m+1})$ $(m\geq 1)$, $^2G_2(3^{2m+1})$
$(m\geq 0).$ Then $\kappa (G) = 2$.
\end{prop}

\begin{remark}
Obviously, it is enough to prove that $\kappa(G_1)=2$ for some
group $G_1$ lying between $G$ and its simply connected cover. In
each specific case the choice of $G_1$ will depend on the
convenience of the proof. In particular, we shall often assume the
Chevalley group under consideration to be simply connected.  We
shall use this observation without any special notice.
\end{remark}

We start with computations for simple groups of Lie type of small
Lie rank defined over the finite fields of size 2, 3, 5. They will
be used in several parts of our proof. The computations were made
for all groups of rank 1 and 2 and also for certain groups of rank
3 and 4 needed for our arguments. The results of MAGMA
computations for groups over ${\mathbb F}_2$ and ${\mathbb F}_3$
are exhibited in Table 1 (over ${\mathbb F}_5$ these groups
contain no nontrivial 2-radical elements). Each entry of Table 1
displays the number of nontrivial 2-radical elements in the
corresponding group (up to conjugacy) and their orders (in
parentheses). Dash means that the corresponding group either is
solvable or does not exist (for this reason the types $A_1$ and
${}^2B_2$ do not appear at all). Asterisks mean that the
corresponding group $G$ is not simple, and computations were made
for the derived subgroup $G'$, which is simple. It is worth
recalling the isomorphisms $B_2(q)\cong C_2(q)$, $B_2(3)\cong
C_2(3)\cong {}^2A_3(2^2)$, $G_2(2)'\cong {}^2A_2(3^2)$,
$A_3(2)\cong A_8$. Note that our computations show that whenever
nontrivial 2-radical elements do exist, there are no nontrivial
3-radical elements.

\begin{table}[!htbp] \label{table1}
\begin{tabular}[h]{|c||c|c||l|}
\hline
{}         &          ${\mathbb F}_2$ & ${\mathbb F}_3$ & Remarks\\
\hline
${}^2A_2$  &  ---                     &      0   &   \\
${}^2G_2$  &  ---                     &      $0^*$   &   (*) Computed for $G'$ \\
$A_2$      &  0                        &     0         & {} \\
$B_2$      &  $0^{**}$                   &     3 (2,3,3)        & (**) Computed for $G'$  \\
$C_2$      &  $0^{**}$                   &     3 (2,3,3)        &  \\
$G_2$      &  $0^{***}$                   &    0     &   (***) Computed for $G'$ \\
${}^2A_3$  &  3 (2,3,3)                &     0      & {}   \\
${}^2A_4$  &  3 (2,3,3)                &     0       & {}  \\
${}^3D_4$  &  0                        &     0       & {}  \\
$A_3$      &  0                        &     0        & {}  \\
$B_3$      &  1 (2)                    &     1 (2)      & {}    \\
$C_3$      &  1 (2)                    &     2 (3,3)        & {}    \\
${}^2F_4$  &  $0^{****}$                        &     ---   & (****) Computed for $G'$     \\
$D_4$      &  0                        &     0       & {}  \\
\hline
\end{tabular}
\caption{2-radical elements in simple groups of small Lie rank}
\end{table}

Before starting the proof of the proposition, we recall the following result from
\cite{Gow} (compare with \cite{EG2}) regarding conjugacy classes of
semisimple elements in Chevalley groups. This fact is essential for our arguments.

\begin{theorem}\label{th:Gow} \cite{Gow} Let $G$ be a finite simple group of Lie type, and
let $g\neq 1$ be a semisimple element in $G$. Let $C$ be a
conjugacy class of $G$ consisting of regular semisimple elements.
Then there exist a regular semisimple $x\in C$ and $z\in G$ such
that $g=[x,z]$.
\end{theorem}

Let us now go over to the proof of Proposition \ref{prop:rank1}.

\begin{proof}
First note that for the groups  $G ={}^2A_2(3^2)$ and  $G
={}^2G_2(3^2)$ the statement of the proposition follows from
calculations presented in Table 1. So we exclude these groups from
consideration in the rest of the proof. We start with several
simple lemmas (recall that $G$ is a {\it finite} group).

\begin{lemma} \label{lem:1-Bruhat}
Let $G = B \cup B\dot w B$ be a group of rank one. Let $1\ne u \in
U$. If $g u g^{-1} \in U$, then $g \in B$.
\end{lemma}

\begin{proof}
Suppose  $g = u_2\dot w u_1$ where $u_1, u_2 \in U$. Then $v = u_2
\dot w u_1 u u_1^{-1} \dot w^{-1} u_2^{-1} \in U$. Hence
$$U \ni u_2^{-1} v u_2 =\dot w u_1 u u_1^{-1} \dot w^{-1} \in
U^-.$$ This contradicts the assumption $u\ne 1$.
\end{proof}

\begin{lemma} \label{lem:oneborel}
Let $G$ be a group of rank one. Then every nontrivial unipotent
element is contained in only one Borel subgroup.
\end{lemma}

\begin{proof}
Suppose $1\ne u \in U \leq B$ and $ u \in B^\prime$, where $B^\prime = x B
x^{-1}, x \notin B$ \cite{St}. Then $u = xvx^{-1}$ for some $v
\in U$. By Lemma \ref{lem:1-Bruhat}, $x \in B$, contradiction.
\end{proof}

\begin{lemma} \label{lem:alternative}
Let $G$ be a group of rank one. Then, up to conjugacy, for every
$g\in G$ we have either $g\in T$, or $g\in U$, or $g=tu$ with $t
\in T, u \in U, tu= ut$, or $g$ is a regular semisimple element
which is not contained in any Borel subgroup.
\end{lemma}

\begin{proof}
Indeed, let $g = su=us$  be the Jordan decomposition of $g$. We
may and shall assume $u\in U$. If $s=1$, then $g=u \in U$, so we
assume further $s \ne 1$. Suppose $u \ne 1 \in U$. Then $s us^{-1}
= u \in U$ and therefore, by Lemma \ref{lem:1-Bruhat}, we have $s
\in B$. Since $s\in B$, $s$ lies in some quasisplit torus. As all
quasisplit tori are conjugate \cite{Ca2}, we have $s'=bsb\minus
\in T$ for some $b\in B$. Thus we get
$$
bgb\minus = bsb\minus bub\minus = s'u'
$$
with $s'\in T$, $u'\in U$. Suppose now $u=1$. We have $g=s$, and
if $s$ lies in a Borel subgroup, then $s$ is conjugate to an
element of $T$, as above. Finally, if $s$ is a semisimple element
which does not belong to any Borel subgroup, then according to
Lemma \ref{lem:1-Bruhat} it does not commute with any unipotent
element, and thus $g=s$ is a regular semisimple element.
\end{proof}

\begin{defn} \label{def:kvadrat}
Let $t \in T$. Define
$$t^{[2]}:= \dot w t^{-1} \dot w^{-1} t.$$
\end{defn}

If $G$ is of the type $A_1$, $^2B_2,$ or $^2G_2$, we have $t^{[2]} =
t^2$. If $G$ is of the type $^2A_2(q^2)$ and $t = \di( \lambda,
\lambda^{-1}\lambda^q, \lambda^{-q})$, we have $t^{[2]} = \di
(\lambda\lambda^q, 1, \lambda^{-1}\lambda^{-q}).$

\begin{lemma} \label{lem:kvadrat}
Let $G$ be a group of rank one, let $g \notin Z(G)$, and let $t$
be a generator of $T$. Suppose $t^{[2]}$ is a regular element.
Then there exists $x \in G$ such that $[g,x]$ is of the form
$\rho^{[2]}$ where $\rho$ is a generator of a quasisplit torus of
$G$.
\end{lemma}

\begin{proof}
We may assume $g = u\dot w$. Put $x = t^{-1}$. Then $$\sigma = [g,
t^{-1}] = u\dot w t^{-1} \dot w^{-1}tt\minus u^{-1} t =
ut^{[2]}t\minus u\minus t$$ which is conjugate to $t^{[2]}v$ for
some $v\in U$.
Since $t^{[2]}$ and, correspondingly, $t^{-[2]}$ are regular
elements, there exists $y$ such that $v=[t^{-[2]},y]$ (see, for
example, \cite{EG2}). Then $yt^{[2]}y^{-1}=
t^{[2]}[t^{-[2]},y]=t^{[2]}v$. Put $\rho=yty^{-1}$. Then $\rho$ is
a generator of a quasisplit torus $T'=yTy^{-1}$ and $w_1=y\dot w
y^{-1}$ is a preimage of the generator of the Weyl group. We have

$$
\begin{array}{ccl}
yt^{[2]}y^{-1} & = & y\dot w t^{-1}\dot w^{-1}ty^{-1}=(y\dot w
y^{-1})(yt^{-1}y^{-1})(y\dot w^{-1}
y^{-1})(yty^{-1}) \\
{} & = & \dot w_1\rho^{-1}w_1^{-1}\rho=\rho^{[2]}.
\end{array}
$$
\end{proof}

\begin{remark}
Explicit calculations with the matrices
$$
\begin{array}{ccl}
t & = & \di(\lambda, \lambda^{-1}), \\
t & = & \di( \lambda, \lambda^{-1}\lambda^q, \lambda^{-q}), \\
t  & = & \di({\lambda},\lambda^{2\theta-1},
\lambda^{-1},\lambda^{1-2\theta}),\\
t & = & \di (\lambda^\theta, \lambda ^{1-\theta}
,\lambda^{2\theta-1}, 1, \lambda^{1-2\theta},\lambda^{\theta
-1},\lambda^{-\theta}),
\end{array}
$$
corresponding, respectively, to the natural representation of
$SL_2(q)$, natural representation of $SU_3(q^2)$, 4-dimensional
representation of the Suzuki group  and 7-dimensional
representation of the Ree group, show that the hypothesis of Lemma
\ref{lem:kvadrat} holds for every group from Proposition
\ref{prop:rank1} except for $A_1(5)$ and ${}^2A_2(3^2)$. These two
groups were considered separately (note that $PSL_2(5) \cong
A_5$).
\end{remark}

\begin{lemma} \label{lem:normalizer}
Let $T'$ be a quasisplit torus in a group $G$ of rank 1, and let $S$
be a subgroup of $T'$ such that $C_G(S)=T'$. Then $N_G(S)=N_G(T')$.
\end{lemma}

\begin{proof}
Let $B^\prime = T^\prime U^\prime$ be a Borel subgroup containing
$T^\prime$, and let $G = B^\prime \cup B^\prime \dot w^\prime
B^\prime$ be the corresponding Bruhat decomposition. Let $g \in
N_G(S)$. Suppose $g = u_1\dot w^\prime u_2$ where $u_1, u_2\in
U^\prime$. Then for every $s \in S$ we have
$$g s g^{-1} = (u_1\dot w^\prime u_2 ) s (u_2^{-1}\dot w^{\prime -1}u_1^{-1}) =
s^\prime \in S \Rightarrow $$$$\Rightarrow B^{\prime -}\ni(\dot
w^\prime s \dot w^{\prime -1}) (\dot w^\prime [s^{-1},u_2] \dot
w^{\prime -1}) = s^\prime [s^{\prime -1},u_1^{-1}]\in B^\prime
\stackrel{B^{\prime -1}\cap B^\prime = T^\prime}{\Rightarrow}$$
$$\Rightarrow [s^{-1},u_2] = 1, [s^{\prime
-1},u_1^{-1}] = 1\stackrel{C_G(S) = T^\prime}{\Rightarrow} u_1 =
u_2 = 1 \Rightarrow g = \dot w^\prime.
$$
Suppose $g \in B^\prime$. Then $g = tu$ for some $t\in T^\prime,
u\in U^\prime$, and for every $s \in S$ we have
$$g s g^{-1} = s t[s^{-1}, u]t^{-1} \in S\Rightarrow [s^{-1}, u] =
1\stackrel{C_G(S) = T^\prime}{\Rightarrow} u = 1.$$ Hence $g \in
N_G(T^\prime)$ and therefore $N_G(S) \leq N_G(T^\prime)$.

Further, using the same arguments as above (put $S = T^\prime$)
one can see that $N_G(T^\prime) = \langle T^\prime ,\,\dot
w^\prime\rangle$. Note that the conjugation with $w^\prime$ is an
automorphism of $T^\prime$ and $T^\prime$ is a cyclic group. Hence
the conjugation with $w^\prime$ is an automorphism of $S$. Thus
$N_G(T^\prime) \leq N_G(S).$
\end{proof}

\begin{lemma} \label{lem:regular}
Suppose the hypothesis of Lemma $\ref{lem:kvadrat}$ holds. Suppose
that for every non-regular $s\in T$ and for every regular $t\in T$
the element $st$ is regular. Then for every $g\notin Z(G)$ there
exist $x, y \in G$ such that the group $H$ generated by $\tau =
[g,x]$ and $\sigma = [g, y]$ is not contained in any Borel
subgroup. Moreover, $\tau \notin N_G(\langle\sigma\rangle)$.
\end{lemma}

\begin{proof}
We shall divide the proof into two cases: 1) $g$ is not a regular
semisimple element; 2) $g$ is a regular semisimple element. Case
1, in turn, will be subdivided into two subcases: 1a)
$\ch(K)\ne2$; 1b) $\ch(K)=2$.

{\it Case 1a)} First suppose $g$ is not a regular semisimple
element. By Lemma \ref{lem:alternative}, we have $g\in B$, $g = s
u $ with $su=us$, $s\in T$ is a non-regular element, and $u \in
U$. Then we can get $1\ne \tau=[g,x] \in U$. Indeed, if $u \ne 1$,
we take $x=s_1\in T$ such that $[u,s_1]\neq 1$. Then
$[g,x]=[g,s_1]=
[us,s_1]=[s,s_1]^u[u,s_1]=[u,s_1]=us_1u^{-1}s^{-1}_1\in U$. If $u
= 1$, then $s \notin Z(G)$, and hence $1\ne [s, v] \in U$ for some
$v \in U$.

Then by Lemma \ref{lem:kvadrat}, we get $\sigma = [g, y] =
\rho^{[2]}$ where $\rho$ is a generator of a quasisplit torus.
Suppose $\left<\tau, \sigma\right>=H \leq B^\prime$ for some Borel
subgroup $B^\prime$. Since $\tau$ is a unipotent element, by Lemma
\ref{lem:oneborel} we have $B^\prime = B$ and therefore
$gyg^{-1}y^{-1} = \sigma\in B$. Consider the element $g^{-1}\sigma
= u^{-1}s^{-1}\sigma$. Since $\sigma  \in B$, we have
$\sigma=s'u'$ where $s'\in T$ is semisimple and $u'\in U$. Since
$\sigma$ is regular, so is $s'$. Then $g^{-1}\sigma =
u^{-1}s^{-1}s'u{'}=s^{-1}u^{-1}s'u'=s^{-1}s'u_1u'=s^{-1}s'u''$ for
some $u{''}\in U$. By the hypothesis of the lemma, $s^{-1}s'$ is a
regular semisimple element. Hence $g^{-1}\sigma = s^{-1}s'u''$ is
a regular semisimple element. Contradiction, since $y g^{-1}
y^{-1} = g^{-1}\sigma$ is not a regular semisimple element.

Let us now prove that $\tau \notin N_G(\langle \sigma \rangle)$.
Assume the contrary. Since $\sigma$ is a regular semisimple
element, we have $C_G(\sigma)=C_G(\left<\sigma\right>)=T'$. Lemma
\ref{lem:normalizer} gives $N_G(\left<\sigma\right>)=N_G(T')$.
Therefore $\tau \in N_G(T^\prime)$.

Hence $\tau^2 \in T^\prime$. Indeed, since $\tau \in
N_G(T^\prime)$, we have $\tau= \dot w^\prime$ where $\dot
w^\prime$ is a preimage of an element of the Weyl group (possibly,
$w^\prime = 1)$ corresponding to $T^\prime$. Thus $\tau^2\in
T^\prime$. But $\tau \in U$. Hence $\tau$ is a unipotent element
of order 2 which contradicts to the assumption $\ch(k) \neq 2$.

{\it Case 1b)} Suppose $g$ is not a regular semisimple element and
$\ch(k)=2$. In this case we may assume $g = \dot w$.

Indeed, let $g = su$ be the Jordan form for $g$. Suppose the order
of $u$ is greater than $2$. On setting $x=t\in T$, we get the
element $[g,x] \in U$ of order greater than $2$. Then, by the
arguments of Case 1a, we have $\tau \notin
N_G(\langle\sigma\rangle)$. Thus the order of $u$ is one or two.
As $\ch(k)=2$, every non-regular element of $T$ lies in the centre
of $G$, and therefore we may assume $s=1$. Hence we may assume $g
= u$ to be an element of order $2$.

As $\ch (k)=2$, in each of the Lie rank one groups, $SL_2(2^m)$,
$SU_3(2^{2m})$, $^2B_2(2^{2m+1})$, all involutions are conjugate,
and we may assume $g = \dot w $.

Therefore we can take $\sigma = [g, t] = [\dot w, t]= t^{[2]},$
and $ \tau = [g,u]=[\dot w, u]=\dot w u \dot w^{-1}
u^{-1}=vu^{-1}$ where $u \in U$ and $1\neq v\in U^-$. Suppose
$\sigma, \tau \in B^\prime$ for some Borel subgroup $B^\prime$.
Then $T \leq B^\prime$ and therefore $B^\prime = B$ or $B^\prime =
B^-$. Contradiction, since  $\tau \notin B$, $\tau \notin B^-$.

Suppose now $\tau = vu^{-1} \in N_G(\langle \sigma \rangle) =
N_G(T)$. This is impossible:
$$vu^{-1}t uv^{-1} = t^{\prime}\in T \Rightarrow
(B\setminus T) \ni u^{-1}tu = v^{-1}t^{\prime} v \in (B^-\setminus
T).$$

{\it Case 2}. Let $g$ be a regular semisimple element. By
\cite{Gow}, we can get $\sigma = [g,y]$ to be a generator of a
quasisplit torus and $\tau = [g,x]$ to be a regular semisimple
element which is not contained in any Borel subgroup.

We have

$\mid T \mid = q-1$ if $G = SL_2(q)$;

$\mid T \mid = q^2-1$ if $G = SU_3(q^2)$ or $G$ is a Suzuki or a
Ree group.

Further,

$(q+ 1)$ divides $\mid G\mid$ if $G = SL_2(q)$, $(q +1, q-1) = 2$
or $1$ (if $q$ is even);

$(q^2-q+1)$ divides $\mid G\mid$ if $G = SU_3(q^2)$, $(q^2-1,
q^2-q+1)$ equals 3 or 1 (indeed, $p$ divides $(q-1)$ implies
$q\equiv 1 \pmod p$, hence $q^2 - q + 1 \equiv 1\pmod p$).
Correspondingly, $ p \mid q +1$ implies $(q^2-q+1) \equiv 3 \pmod
p$;

$(q^4+1)$ divides $\mid G\mid$ if $G$ is a Suzuki group, $q^2 =
2^{2m+1}$, $(q^2-1, q^4+1) = 1$;

$(q^4-q^2+ 1)$ divides $\mid G\mid$ if $G$ is a Ree group, $q^2 =
3^{2m+1}$, $(q^4-q^2+1, q^2-1) = 1$.

Let now $G = SL_2(q)$. Then the maximal nonsplit torus is a cyclic
group of order $q+1$. By \cite{Gow}, we can take $\tau = [g, y]$
to be a generator of such a group. Then the order of $\tau^2$ is
equal to $q+1> 2$ if $q =2^m$ or  $(q+1)/2 > 2$ (note that $q
>3$). Hence $\tau \notin N_G(\langle \sigma \rangle) = N_G(T)$
(because $\tau^2 \notin T$). Also $\tau$ does not belong to a
Borel subgroup.

Let $G = SU_3(q^2)$. Suppose that 3 divides $q^2-q +1$. Then
$q\equiv -1 \pmod 3$, hence $q \equiv 2, 5, 8 \pmod 9$ and,
therefore, $ 9$ does not divide $ q^2 -q + 1$. Then there exists a
prime $p \ne 2, 3$, $p \mid q^2-q +1 $. By \cite{Gow}, we can
obtain an element of order $p$ of the form $\tau = [g, y]$. Then $
\tau \notin N_G(\langle \sigma \rangle)$, and $\tau$ does not
belong to a Borel subgroup.

If $G$ is of Suzuki or Ree type,  take $p\mid q^4 +1$ or $p\mid
q^4-q^2+1$, respectively, and proceed as above.

Thus, in all the cases $\tau \notin N_G(\langle \sigma \rangle)$.
\end{proof}

\begin{remark}
The hypotheses of  Lemma \ref{lem:regular} hold for every group from
Proposition \ref{prop:rank1}. This can also be checked by explicit
calculations with diagonal matrices (see \cite{Ca2} and \cite{KLM}).
\end{remark}

\begin{lemma} \label{lem:nonsolv1}
There exist $\tau = [g,x]$ and $\sigma = [g, y]$ such that the
subgroup $H=\left<\sigma ,\tau\right>$ is not solvable.
\end{lemma}

We choose $\tau = [g,x]$ and $\sigma = [g, y]$ as in the previous
lemma.

It is enough to show that $H$ does not contain abelian normal
subgroups. Let $A$ be a maximal abelian normal subgroup of $H$. We
want to check that $A$ is a reductive group. Suppose $p = \char
(K) $ divides the order of $A$. Then the Sylow $p$-subgroup of $A$
is normalized by $H$. By Lemma  \ref{lem:alternative}, $H \leq
B^\prime$ for some Borel subgroup $B^\prime$. This is impossible
in view of Lemma  \ref{lem:regular}. Hence the order of $A$ is not
divisible by $p$, and $A$ is a reductive group.

Let us now view $H$ as a subgroup of $GL(V)$ where $V$ is a finite
dimensional vector space over an algebraically closed field and
$\dim V = 3$ (if $G = PSL_2(q), q \ne 2^n$), $\dim V = 2$ (if $G =
SL_2(2^m)$), $\dim V = 8$ (if $G= PSU_3(q^2),$ $\dim V = 4$ (case
$^2B_2$), or $\dim V =7$ (case $^2G_2$). Then $A$ is
diagonalizable in $GL(V)$ and not all irreducible components of
the $A$-module $V$ are isomorphic (if $A \ne Z(H) $). Thus there
exists a nontrivial homomorphism $\varphi \colon H \rightarrow
S_k$, $ k \leq 3, 2, 8, 4, 7$ which corresponds to permutations of
isotypical components (otherwise, $A \leq Z(H)$).

\medskip

{\it Case 1}. { \it Let $G = PSL_2(q), q \ne 2^m$.} For $q\leq 25$
the statement of the lemma is checked by explicit computer
calculations with MAGMA. Let now $q>25$. Recall  that $\sigma
=t^2$ or $\sigma = t$ for $\langle t \rangle=T'$, where $T'$ is a
split torus in $G$. Since the order of $T'$ is $\geq (q-1)/2$, the
order of $\sigma$ is  $\geq (q-1)/4$ . Since $\varphi(\sigma)$
lies in $S_3$, we have $\varphi(\sigma^n)=1$ for some $n\leq 3$.
Thus $ord \ \sigma^n\geq (q-1)/12>2$. Hence $C_G(\sigma^n)=T'$
because $\sigma^n$ is a regular semisimple element of $T^\prime$.
Since $\varphi (\sigma^n)=1$, we have $\sigma^n\in C_H(A)$.

\begin{sublemma} \label{sublemma}
i) With the above notation, suppose there exists $h\in H$ such that
\item
1. $h\in C_H(A)$;
2. $h\in T'$;
3. $C_G(h)=T'$.
Then $A\subseteq T'$.

ii) If, in addition, there exists $a\in A$ such that $C_G(a)=T'$,
then $N_G(\left<h\right>)=N_G(A)=N_G(T')$.
\end{sublemma}

\begin{proof}
The first assertion of the sublemma is obvious: if  $h\in C_H(A)$,
then $a\in C_G(h)$ for any $a\in A$. The second assertion follows
from Lemma \ref{lem:normalizer} applied to $S=A$ and
$S=\left<h\right>$.
\end{proof}

On setting $h=\sigma^n$, we conclude that $A\subseteq T'$.

Suppose there exists $a$ generating $A$ such that $C_G(a)=T'$.
Then by the above sublemma we have $N_G(\left< \sigma^n\right>
)=N_G(A)=N_G(T').$ On the other hand, we have
$N_G(T')=N_G(\left<\sigma\right>)$. (Indeed, the inclusion
$N_G(\left<\sigma\right>)\subseteq N_G(\left<\sigma^n\right>)$ is
obvious, and the inclusion $N_G(T')\subseteq
N_G(\left<\sigma\right>)$ follows from the fact that in the groups
of Lie rank 1 the generator $w$ of the Weyl group normalizes $t\in
T$ and hence $\sigma$.) Thus we conclude that
$N_G(\left<\sigma\right>)=N_G(A)\supseteq H$, which contradicts
the choice of $\tau$.

Suppose now there is no $a\in A$ such that $C_G(a)=T'$. Then
$A=\left<a\right>$ is a cyclic subgroup of order 2 (all other
elements of $T'$ are regular). Since the order of $a$ equals 2, we
have  $N_G(A)=C_G(A)$. On the other hand,
$C_G(A)=N_G(T')=N_G(\left<\sigma\right>)$. Again we get
a contradiction since $\tau$ belongs to $H\subseteq N_G(A)$ but does
not belong to $N_G(\langle \sigma\rangle )$.

\medskip

{\it Case 2}. Let $G = SL_2(2^m), m > 1$. In this case $G$ has no
centre, any element of $T'=\left<t\right>$ is regular, the order
of $t$ equals $2^m-1$. Hence the order of $\sigma^2$ equals $2^m-1
> 1$. Therefore we can use the same argument as in the preceding
case.

\medskip

{\it Case 3}. Let $G = PSU_3(q^2), q > 3.$ 
In this case the semisimple element $\sigma= t$ or $t^{[2]}$. The
order of the image of $\sigma $ in $PSU_3(q^2)$ is $\geq q-1$
(recall that the centre of $SU_3(q^2)$ is nontrivial if and only
if $q+1=3k$ for some $k$). Note that $\sigma^n$ is a nonregular
nontrivial element if and only if $\sigma^n=\di(-1,1,-1)$. Hence
if $n \leq 8$ and $q > 17$, the order of $\sigma^n \geq (q-1)/8 >
2$ and therefore the image of $\sigma^m$ in $PSU_3(q^2)$ is a
regular element. Thus we may use the same arguments as in the
previous case. Explicit computer calculations with MAGMA prove the
statement for the remaining cases $q\le 17$.

\medskip

{\it Case 4.} Let $G$ be a Suzuki or a Ree group. Every nontrivial
element of $T$ is regular if $G$ is a Suzuki group \cite{Ca2}, and
every element of $T$ of order greater than two is regular if $G$
is a Ree group \cite{KLM}.  Note that if $G$ is a Ree group, then
the order of a maximal torus $T^\prime$ is equal to $3^{2m+1}-1$.
Hence $2 \mid T^\prime\mid, 4 \nmid \mid T^\prime \mid $. The
element $\sigma$ is a generator or the square of a generator of
$T^\prime$. In particular, $\sigma$ is not an involution. So if
$n$ is less than the order of $\sigma^2$, then $\sigma^n$ is a
regular element of a maximal quasisplit torus.

Consider the permutation $\varphi(\sigma )\in S_k$. First suppose
$\varphi(\sigma )=1$. Arguing as in Case 1, we arrive at a
contradiction with the choice of $\tau$ whenever we can choose
$a\in A$ such that $C_G(a)=A$. This is always possible except for
the case where $G$ is a Ree group and $A$ is generated by the
(unique up to conjugacy) involution $a$ of $G$. But in this latter
case we have $N_G(A)=C_G(a)=\mathbb Z/2\times PSL_2(3^{2m+1})$
\cite[Th.~3.33(iv)]{Gor2}. Hence $H\subseteq PSL_2(3^{2m+1})$, and
we are reduced to Case 1.

Thus we may assume $\varphi(\sigma )\ne 1$. Then the same argument
as above with $\sigma^n$ replacing $\sigma$ shows that
$\varphi(\sigma^n) \ne 1$ for every $n < ord \sigma^2$. This means
that the restriction of $\varphi$ to $\left< \sigma^2\right>$ is
faithful. But this is impossible since $\varphi(\sigma^2) \in S_4$
for the Suzuki groups and the order of $\varphi(\sigma)$ must be
less than or equal to 4. However in this case $ord
\varphi(\sigma^2)=ord \sigma^2 = 2^{2m+1}-1 > 4$. The same
situation takes place for the Ree groups: $ord
\varphi(\sigma^2)=ord \sigma^2 = (3^{2m +1}-1)/2 > 12$, and
therefore $\varphi(\sigma^2)$ cannot belong to $S_7$.

Thus in the Suzuki and Ree groups there are no nontrivial abelian
normal subgroups in $H$, and hence $H$ is not solvable.

\medskip

Lemma \ref{lem:nonsolv1} (and hence Proposition \ref{prop:rank1})
are proved.
\end{proof}

To use Proposition \ref{prop:rank1} as induction base, we have to
extend it from the {\it simple} groups of rank 1 to some {\it
reductive} groups of semisimple rank 1, namely to the case where
$G$ is an extension of a simple group by a diagonal automorphism,
because such groups appear as Levi factors of parabolic subgroups
of simple groups of higher ranks; see Lemma \ref{lem:Levi} below.

\begin{prop} \label{prop:rank1-red}
Let $G$ be one of the groups from the list of Proposition
$\ref{prop:rank1}$, let $h$ be a diagonal automorphism of $G$, and
let $L$ denote the corresponding extension. Then $\kappa (L)=2$
with the sole exception $L=PGL_2(5)$ for which we have $\kappa
(L)=3$.
\end{prop}

\begin{proof} The exceptional case is not a surprise in light of
the isomorphism $PGL_2(5)\cong S_5$. In order not to overload the
reader with technicalities, we shall only sketch the proof. The
key point is the following generalization of Gow's Theorem
\ref{th:Gow}:

\begin{theorem} \label{th:Gow-red} \cite{Gord}
Let $L = \tilde{L}^F$ be the group of fixed points of the
Frobenius map $F$ acting on a connected reductive group
$\tilde{L}$ defined over $\bar{\mathbb F}_p$ such that the derived
subgroup $\tilde{L}^\prime$ is a simple algebraic group, and
denote $L'=(\tilde{L'})^F$. Suppose that $([L:L'],p)=1$. Let
$\gamma, g\in L$ be semisimple elements such that $\gamma$ is
regular, $g \in L'$, $g \notin Z(L)$. Then there exist $g' \in
C_g$ and $x \in L$ such that $g' = [\gamma, x].$ \qed
\end{theorem}

Note that the hypotheses of the theorem hold if $L'$ is generated
by the root subgroups of $L$. The proof goes along the same lines
as in \cite{Gow}.

Furthermore, Lemmas \ref{lem:1-Bruhat}--\ref{lem:alternative} hold
for all reductive groups of semisimple rank 1, and Lemmas
\ref{lem:kvadrat}, \ref{lem:regular}, \ref{lem:nonsolv1} admit
appropriate modifications. For example, the statement of Lemma
\ref{lem:kvadrat} should be modified as follows:

\begin{lemma} \label{lem:kvadrat-red}
Let $L=HG$ be a reductive group of semisimple rank one, where $G$
is the derived subgroup of $L$ and $H$ is generated by a diagonal
automorphism $h$ of $G$. Let $g\in L$, $g \notin Z(G)$, and let
$t$ be a generator of $T$. Suppose $t^{[2]}$ is a regular element.
Then there exists $x \in G$ such that $[g,x]$ is of the form
$\rho^{[2]}$ where $\rho$ is a generator of a quasisplit torus of
$G$. \qed
\end{lemma}

Note that to prove the modified Lemma \ref{lem:nonsolv1},
additional MAGMA computations are needed to treat the reductive
groups $PGL_2(q)$ ($q\le 25$) and $PGU_3(q)$ ($q\le 17$).
\end{proof}



\section{Groups of Lie type of arbitrary rank over fields of large characteristic}\label
{sec:lie_gen}

\begin{theorem} \label{th:bigrank}
Let $G$ be a Chevalley group of rank  $ > 1$ over field $\ch (K)
\ne 2, K \ne {\mathbb F}_3$. Then $\kappa(G) = 2$.
\end{theorem}

\begin{proof}
We need several lemmas (most of whose statements are independent
of the characteristic of the ground field).

\begin{lemma} \label{lem:positive}
Let $\Pi =\{\alpha_1,\dots, \alpha_r\}$, $r \geq  2, $ be a basis
of an irreducible root system $R\neq A_2$, where the numbering of
the simple roots is as in \cite{Bou} in the case $R \ne E_r$, and
$\alpha_2$ and $\alpha_3$ are interchanged in the case $R = E_r$.
Denote by $w_c = w_{\alpha_1}\cdots w_{\alpha_{r}} w_{\alpha_2}$
the corresponding Coxeter element. Then $w_c(\alpha_1) > 0,\
w_c(\alpha_1) \notin \Pi$ and $w^{-1}_c(\alpha_2) > 0,\
w^{-1}_c(\alpha_2) \notin \Pi$.
\end{lemma}

\begin{proof}

Let $r=2.$ We proceed case by case.

1. $R = B_2.$
 We have $\alpha_1 = \epsilon_1 -
\epsilon_2,\ \alpha_2 = \epsilon_2 ,$
 and
$$w_c(\alpha_1) = \epsilon_1 + \epsilon_2=\alpha_1+2\alpha_2,\quad w_c^{-1}(\alpha_2) =
\epsilon_1=\alpha_1+\alpha_2.$$

2. $R = C_2.$
 We have $\alpha_1 = \epsilon_1 -
\epsilon_2,\  \alpha_2 = 2\epsilon_2 ,$
 and
$$w_c(\alpha_1) = \epsilon_1 + \epsilon_2=\alpha_1+\alpha_2,\quad w_c^{-1}(\alpha_2) =
2\epsilon_1=2\alpha_1+\alpha_2.$$

3. $R = G_2$. Then $\alpha_1 =\epsilon_1 - \epsilon_2, \ \alpha_2 =
-2\epsilon_1 + \epsilon_2+\epsilon_3$. We have
$$w_c(\alpha_1) = \epsilon_3 -\epsilon_2 = 2\alpha_1 + \alpha_2
,\,\,\, w_c^{-1}(\alpha_2) = 2\epsilon_3 -\epsilon_2 -\epsilon_2 =
3\alpha_1 + 2\alpha_2.$$

Let $r > 3. $ Note that our numbering of roots gives
$\langle\alpha_1,\alpha_2\rangle = A_2$. Therefore
\begin{equation}
w_{\alpha_1}(\alpha_2) = \alpha_1+\alpha_2,
\,\,\,w_{\alpha_2}(\alpha_1) = \alpha_1+\alpha_2,\label{eq:21}
\end{equation}
\begin{equation}
w_{\alpha_1}(\alpha_1 +\alpha_2) =\alpha_2, \ \
w_{\alpha_2}(\alpha_1 +\alpha_2) =\alpha_1.\label{eq:22}
\end{equation}
Put $\omega = w_{\alpha_{3}}\cdots w_{\alpha_{r}}$. Since $\omega$
has no factors $w_{\alpha_{1,2}}$,  we have
\begin{equation}
\omega^{\pm 1}(\alpha_{1,2})
> 0.\label{eq:23}
\end{equation}
Moreover,
\begin{equation}
\omega^{\pm 1}(\alpha_1) = \alpha_1,\ \ \omega^{\pm 1}(\alpha_2 )
\notin \langle\alpha_1,\alpha_2\rangle.\label{eq:24}
\end{equation}

From (\ref{eq:21})--(\ref{eq:24}) we get
\begin{equation}
\omega^{\pm 1}(\alpha_1 +\alpha_2) = \alpha_{1}+ \alpha_2 +...\ne
\alpha_1 +\alpha_2,\ \omega^{\pm 1}(\alpha_1 +\alpha_2)
> 0.\label{eq:25}
\end{equation}

From (\ref{eq:25}) we get $$\begin{aligned}0 <  w_c(\alpha_1) =
w_{\alpha_1}\omega(\alpha_1 +\alpha_2)\notin \Pi,\\0 <
w^{-1}_c(\alpha_2) = w_{\alpha_2}\omega^{-1}(\alpha_1
+\alpha_2)\notin \Pi.
\end{aligned}$$
\end{proof}

\begin{lemma} \label{lem:lefttwo}
Let $g  =  u^{-1}\dot w^{-1}_c  ,$ where $w_c$ is the Coxeter
element from the previous lemma and $u \in U$. Then there exists
$x\in G$ such that $[g,x] = u_{\alpha_1}u_{\alpha_2} u^\prime$,
where $u_{\alpha_1} \ne 1, u_{\alpha_2} \ne 1$ are the
corresponding root subgroup elements and $u^\prime \in U$ does not
contain root subgroups factors of type $u_{\alpha_1}, u
_{\alpha_2}$. Moreover, every $u_{\alpha_1}\in U_{\alpha_1}$ can
be obtained in such a way.
\end{lemma}

\begin{proof}
Let $R = A_2$. Put $1\ne x = u^\prime_{\alpha_2}\in U_{\alpha_2}$.
Then $\dot w^{-1}_cu^\prime_{\alpha_2}\dot w_c =
u^\prime_{\alpha_1} \in U_{\alpha_1}$ and
$$[g,x] = u^{-1}(\dot w^{-1}_c u^\prime_{\alpha_2}\dot w_c)
uu_{\alpha_2}^{\prime -1}  = (u^{-1} u^\prime_{\alpha_1}
u)u_{\alpha_2}^{\prime-1} = u_{\alpha_1}u_{\alpha_2} u^\prime$$
where $u_{\alpha_1} = u_{\alpha_1}^\prime, u_{\alpha_2} =
u_{\alpha_2}^{\prime-1}$, and $u^\prime = u_{\alpha_2}^\prime
[u_{\alpha_1}^{\prime-1}, u^{-1}] u_{\alpha_2}^{\prime -1}$ does
not contain factors from $U_{\alpha_1}$, $U_{\alpha_2}$.

On varying $x=u'_{\alpha^\prime_2}$, we can get an arbitrary
$u_{\alpha_1}$.

Let now $R \ne A_2$. We use Lemma \ref{lem:positive}. Put  $x
=u^\prime_{\alpha_2}u^\prime_{\beta}$ where $\beta = w_c(\alpha_1)$.
Then $\dot w^{-1}_cu^\prime_{\alpha_2}\dot w_c =
u^\prime_{\gamma},\gamma
> 0, \gamma \notin \Pi,\,\,\,\dot w^{-1}_cu^\prime_{\beta}\dot w_c = u^\prime_{\alpha_1} \in
U_{\alpha_1}$, and
$$[g,x]= u^{-1}(\dot w^{-1}_c u^\prime_{\alpha_2}u^\prime_\beta\dot w_c)
uu_\beta^{\prime-1}u_{\alpha_2}^{\prime -1} =
 (u^{-1}u^\prime_\gamma
u^\prime_{\alpha_1} u)u_\beta^{\prime-1}u_{\alpha_2}^{\prime-1} =
u_{\alpha_1}u_{\alpha_2} u^\prime ,$$ with $u'$ as required.
\end{proof}

\begin{lemma} \label{lem:leftone}
Let $g  =  u^{-1}\dot w^{-1}_c  ,$ where $w_c$ is the Coxeter
element from Lemma $\ref{lem:positive}$. Then there exists $y\in
G$ such that $[g,y] = u_{-\alpha_1} u^\prime$ where $u_{-\alpha_1}
\in U_{-\alpha_1},\,\,\,u^\prime \in U$. Moreover, every
$u_{-\alpha_1}\in U_{-\alpha_1}$ can be obtained in such a way.
\end{lemma}

\begin{proof}
Put  $y = u^{-1}_{-\alpha_1}$. We have $\dot
w^{-1}_cu^{-1}_{-\alpha_1}\dot w_c = u_\beta,$ $\beta
>0$, and $\beta\neq \alpha_1$ (this follows from the
definition of $w_c$). Then
$$
\begin{array}{ccl}
[g,y] & = & u^{-1}\dot w^{-1}_c u^{-1}_{-\alpha_1}\dot w_c
uu_{-\alpha_1} = u^{-1}(\dot w^{-1}_cu^{-1}_{-\alpha_1}\dot w_c)
uu_{-\alpha_1}\\
{} & = & u^{-1}u^{-1}_{\beta} uu_{-\alpha_1} =
u_{-\alpha_1}(u_{-\alpha_1}^{-1}u^{-1}u^{-1}_{\beta}
uu_{-\alpha_1}) = u_{-\alpha_1} u^\prime.
\end{array}
$$
The last equality follows from the fact that $u\minus
u\minus_\beta u$ belongs to the unipotent radical of the minimal
parabolic subgroup corresponding to the root $\alpha_1$.
\end{proof}

\begin{lemma} \label{lem:DN}
Let $G$ be a quasisimple Chevalley group of rank one $\ne A_1(2^m).$
Then there exist $u_1 \in U^-, u_2 \in U^+$ such that $\langle u_1,
u_2\rangle$ is not solvable.
\end{lemma}

\begin{proof}

The proof immediately follows from Dickson's lemma (see
\cite[Theorem 2.8.4]{Gor2} and \cite{Nu}), where there are exhibited
explicit pairs of unipotent elements $u_1 \in U^-, u_2 \in U^+$ such
that the subgroup $\langle u_1, u_2\rangle $ is not solvable.
\end{proof}

The following lemma allows us to reduce to groups of small
semisimple rank.

\begin{lemma} \label{lem:Levi}
Let $P = LV$ be a parabolic subgroup of a Chevalley group $G$
where $L$ is a Levi factor and $V$ is the unipotent radical of
$P$. Further, let $x_1, \ldots, x_s,g \in P$, and let $\bar{x}_1,
\dots, \bar{x}_s, \bar{g}$ be their images in $L/Z(L)$ with
respect to the natural homomorphism $P \rightarrow L\rightarrow
L/Z(L)$. If the group $\langle [\bar{g},
\bar{x}_1],\dots,[\bar{g},\bar{x}_s]\rangle$ is not solvable, then
the group $\langle [g,x_1],\dots,[g,x_s]\rangle$ is not solvable
too.
\end{lemma}

\begin{proof}
Obvious.
\end{proof}

\smallskip

Now we are able to finish the proof of the theorem. Let $X \subset
\Pi$, and let $ X = X_1\cup\ldots\cup X_l$ be the decomposition of
$X$ into a disjoint union of subsets $X_i$ generating irreducible
subsystems of $R$. Put
\begin{equation}
w_{X_i} = \prod_{\alpha \in X_i} w_{\alpha}\label{eq:26}
\end{equation}
where the product is taken in any order. Set $$w_X = \prod_{i}
w_{X_i}.$$ (If $X=\emptyset$, we set $w_X=1$.) Then $w_X$ is a
generalized Coxeter element (see \cite{GS}) corresponding to $X$.
Denote $\dot W_X= \langle \dot w_\alpha, \alpha\in \langle
X\rangle\rangle$, where $\langle X\rangle$ stands for the root
system generated by $X$. Since $\ch(K)\ne 2$, $G$ is not of type
${}^2F_4$, and we can use the following

\begin{prop} \label{prop:cell} \cite[Prop.~6]{GS}
Suppose $G$ is not of type ${}^2F_4$, and  let $g \in G\setminus
Z(G)$. Then the conjugacy class of $g$ intersects a generalized
Coxeter cell $B \dot w_X B$ for some $X$.
\end{prop}

\begin{remark} \label{rem:st}
For $G={}^2F_4$ it is not known whether the above statement is
true or not.
\end{remark}

Thus we may assume
$$g = u \dot w_X,\,\,\,u \in U.$$

To finish the proof of Theorem \ref{th:bigrank}, we now consider
three separate cases. (Note that if $X\neq \emptyset$, we have
$\mid X_i \mid \neq \emptyset$ for every $i$.)

{\it Case 1.} Suppose $X = \emptyset$. Then $g = hu,$ $u\in U,$ $h
\in T$. We may assume $u \ne 1$ (otherwise we can conjugate $g$
with an appropriate element from $U$).  Conjugating $g$ with an
appropriate element $\dot w$ we can get an element $g^\prime =
 h^\prime u^\prime$ in the conjugacy class of $g$ such that
$u^\prime \in U$ and among root factors of $u^\prime$ there is a
simple root subgroup factor $u_{\alpha}$.

Indeed, let $$u = \prod_{\alpha \in M\subset
R^+}u_{\alpha},\,\,\,u_{\alpha}\ne 1.$$ Let $k =
\min\{ht(\alpha)\,\,\,\mid\,\,\,\alpha\in M\}$. Then there exists an
element $w \in W$ such that $0 <
\min\{ht(\alpha)\,\,\,\mid\,\,\,\alpha\in w(M)\} < k$. Thus we can
get $\min\{ht(\alpha)\,\,\,\mid\,\,\,\alpha\in w(M)\} = 1$ for an
appropriate $w\in W$.

Write $g=hu_{\alpha}u''.$

Put $P = T\langle U_{\pm \alpha}\rangle U = B\langle\dot
w_{\alpha}\rangle B$. Now in the parabolic subgroup $P$ we can
take the Levi factor $L_{\alpha}$ of rank 1 corresponding to the
root $\alpha$. Denote by $\bar g$ the image of $g$ in
$\overline{L}_{\alpha}=L_{\alpha}/Z(L_{\alpha})$. First suppose
that $\bar g$ lies in the derived subgroup
$G_{\alpha}=L^{\prime}_{\alpha}$ of $L_{\alpha}$. Since
$G_{\alpha}$ is a simple group of rank 1, we can apply Lemma
\ref{lem:Levi} and Proposition \ref{prop:rank1} and get the
result. If $\bar g$ does not lie in $G_{\alpha}$, we can use
Proposition \ref{prop:rank1-red} instead of Proposition
\ref{prop:rank1} except for the case $p=5$. In this latter case we
have $G_{\alpha}\cong PGL_2(5)$, and we have a problem only when
$\bar g$ is a (unique up to conjugacy) 2-radical element of
$PGL_2(5)$. One can show that in fact this case cannot occur.
Indeed, this 2-radical element is an involution which can be
represented by the matrix $\left( \begin{matrix} 0 & 1\\3 & 0
\end{matrix}\right) \in PGL_2(5)$. One can easily show that this matrix is
not triangulizable which contradicts the form of $g$.

\bigskip

{\it Case 2.} Suppose $\mid X_i\mid > 1$ for every $i$. Put $P =
B\dot W_X B$. Consider the group $L_i = T \langle U_{\pm
\alpha}\,\,\mid\,\,\,\alpha \in \langle X_i\rangle \rangle$. This
is a subgroup of a Levi factor $L = T\langle U_{\pm
\alpha}\,\,\mid\,\,\,\alpha \in \langle X\rangle \rangle$ of $P$.
Let $g_i = u_iw_{X_i}$ be $i^{th}$ component of $g$. We may assume
that the order of simple reflections in (\ref{eq:26}) corresponds
to the order in Lemmas \ref{lem:positive}--\ref{lem:leftone}. Then
by Lemmas \ref{lem:lefttwo}--\ref{lem:leftone}, given any
$u_{\pm\alpha_i}\in U_{\pm\alpha_i}$, we have
$$[g_i ,x] = u_{\alpha_{i_1}}
u^{\prime\prime}, [g_i, y] = u_{-\alpha_{i_1}}u^\prime$$ for some
$x,y\in G$ and $u', u''$ as specified there. It remains to use
Lemmas \ref{lem:DN}--\ref{lem:Levi}; note that Lemmas
\ref{lem:lefttwo}--\ref{lem:Levi} hold for the reductive (not
necessarily simple) groups appearing as Levi factors of parabolic
subgroups.
\end{proof}

\bigskip

{\it Case 3.} Suppose $\mid X_i\mid = 1$ for some $i$. Let $P =
B\dot W_{X}B$. Then  there exists a simple component $L_i$ of a
Levi factor of $P$ which is of rank one. If $p > 5$, we can use
Lemma \ref{lem:Levi} and Proposition \ref{prop:rank1-red}. To
treat the exceptional case $p=5$, we reduce to consideration of
the reductive groups of ranks 2, 3, 4 over ${\mathbb F}_5$ (see
Case 4 at the end of the next section for a more detailed
argument). MAGMA computations show that these groups contain no
nontrivial 2-radical elements, and we are done.

\medskip

Theorem \ref{th:bigrank} is proved.

\section{Groups of Lie type over fields of small characteristic}\label{sec:char_2}

\begin{prop} \label{prop:bigrank}
Let $G$ be a nonsolvable Chevalley group over a field $K$ where
either $\ch (K) = 2$ or $K={\mathbb F}_3$. Then $\kappa (G) \le
3$.
\end{prop}

\begin{proof}

Throughout this section we assume $G\ne {}^2F_4(q^2)$ leaving this
case for separate consideration in the next section.

We have to prove that for every $g\notin Z(G)$ one can find $x_1,
x_2, x_3 \in G$ such that the group $F = \langle [g,x_1], [g,
x_2], [g, x_3]\rangle$ is not solvable.

By Proposition \ref{prop:rank1}, for any nonsolvable rank 1 group
$G$ we have $\kappa (G) = 2$. Thus we may and shall assume that
$\rank \ G > 1$.

First suppose $\ch (K) = 2$, $\mid K \mid >2$. We use the same
case-by-case subdivision as in the proof of Theorem
\ref{th:bigrank} above. Cases 1 and 3 are treated in exactly the
same way (two commutators are enough). Suppose that we are in the
conditions of Case 2, i.e. $\mid X_i\mid > 1$ for every $i$.
Arguing as in the proof of Theorem \ref{th:bigrank}, we reduce to
the case of $L_{\alpha}/Z(L_{\alpha})$, where $L_{\alpha}$ is a
Levi factor of semisimple rank 1. If $L_{\alpha}$ is not of type
$A_1(2^m) (m>1)$, we can use the same arguments as in Lemma
\ref{lem:DN} (once again, two commutators are enough). So we may
and shall assume $G$ of type $A_1(2^m) (m>1)$.

Arguing as in the proofs of Lemmas
\ref{lem:lefttwo}--\ref{lem:leftone}, we conclude that there exist
$x_1, x_2, x_3 \in G$ such that $[g,x_1] = v \in U^-$, $[g, x_2]
=u^\prime \in U,$ $[g, x_3]= u \in U,$ $ u\notin \langle
u^\prime\rangle,$ where $v$, $u'$, $u$ are arbitrary given
elements. Moreover, according to \cite{EG1}, \cite{CEG}, we can
arrange our choice so that to make $ s = v u'$ a generator of a
maximal split torus of $G$. Finally, note that $u$ is a regular
unipotent element (as all unipotent elements in $SL_2(2^m)$).

Put $\sigma = s, \tau = u$.
Since $v, u'$ are  involutions, $u^\prime$ belongs to
$N_G(\langle \sigma\rangle)$. Indeed, we have $u^\prime
\sigma{u^\prime}^{-1}=u^\prime \sigma {u^\prime}=u^\prime s
u^\prime=u^\prime v u^\prime {u^\prime}=u^\prime v=(v
u^\prime)^{-1}=\sigma^{-1}$. Then $\tau=u$ does not belong to
$N_G(\langle \sigma\rangle)$ (otherwise we would have $u,\, u^\prime
\in N_G(\langle \sigma\rangle)$ and $\mid \langle u, u^\prime
\rangle\mid\  =4$, contradiction to $\mid N_G(\langle
\sigma\rangle)\mid\ =2(2^m-1)$).
Further, $u$ and $ vu^\prime$ cannot be in the same parabolic subgroup
($u$ can belong only to $B$ (Lemma \ref{lem:oneborel}), but
$vu^\prime \notin B$). Now we can repeat the arguments used in the
proof for rank one groups over a field of odd characteristic (see
Lemmas \ref{lem:regular} and \ref{lem:nonsolv1}).

\bigskip

Let now $\mid K \mid = 2$ or $\mid K\mid = 3$.

\medskip

{\it Case 1. $X = \emptyset$}. Then $g = uh,$ $u\in U,$ $h \in T$.
We may assume $u \ne 1$ (otherwise we can conjugate $g$ with an
appropriate element from $U$). Conjugating $g$ with an appropriate
element $\dot w$, we can get $g^\prime = u^\prime h^\prime$ in the
conjugacy class of $g$ such that $u^\prime \in U$ and among root
factors of $u^\prime$ there is a nontrivial simple root subgroup
factor $u_{\alpha_i}$ (see Case 1 in the end of the proof of
Theorem \ref{th:bigrank}). Let $\alpha_j$ be any root adjacent (in
the Dynkin diagram) to $\alpha_i$. Then we can reduce to the case
of a Levi factor of semisimple rank 2, as above. For all groups of
rank 2 over ${\mathbb F}_2$ and ${\mathbb F}_3$ we use explicit
MAGMA computations. Table 1 contains the needed data for simple
groups. As to reductive groups of semisimple rank two, 2-radical
elements can only appear in groups of type $C_2$ (isomorphic to
$B_2$) and in $PGU_4(2^2)$ (isomorphic to $PSp_4(3)$); in all
these cases there are no nontrivial 3-radical elements.

\medskip

{\it Case 2. $\mid X_i\mid = 2$ for some $i$}. Let $P$ be a
parabolic subgroup, $P = B\dot W_{X}B$. Then  there exists a
simple component $L_i$ of a Levi factor of $P$ which is of
semisimple rank two. Then we can use Lemma \ref{lem:Levi} and
explicit MAGMA computations for the groups of rank two, as in the
previous case.


\medskip

{\it Case 3. $\mid X_i\mid > 2$ for some $i$}.

In this case, the arguments based on the use of Lemmas
\ref{lem:lefttwo}--\ref{lem:leftone} are not enough. Instead we
shall use the following more subtle version of Lemma
\ref{lem:leftone}. Note that as Lemma \ref{lem:leftone}, this
lemma holds for the reductive, not necessarily simple, groups
appearing as Levi factors of parabolic subgroups.

\begin{lemma} \label{lem:leftonemod}
Let $g  =  u^{-1}\dot w^{-1}_c  ,$ where $w_c$ is the Coxeter
element from Lemma $\ref{lem:positive}$. Then:

1) there exists $y\in G$ such that $[g,y] = u_{-\alpha_1}
u^\prime$ where $u_{-\alpha_1}$ is any prescribed element from $
U_{-\alpha_1}$ and $\,\,\,u^\prime \in U$;

2) there exists $z \in G$ such that $[g,z] = f u^{\prime\prime}$
where $f \in \langle U_{\alpha_2}, U_{-\alpha_2}\rangle, f \notin
B$ and $u^{\prime\prime} \in U$.
\end{lemma}

\begin{proof}
1) See Lemma \ref{lem:leftone}.

2) Recall that $w_c = w_{\alpha_1}\cdots w_{\alpha_{r}}
w_{\alpha_2} = \omega w_{\alpha_2}$. Since $\omega$ does not
contain the factor $w_{\alpha_2}$, we have $\omega (\alpha_2) =
\gamma
> 0$ and $w_c^{-1}(\gamma) = w_{\alpha_2}\omega^{-1}(\gamma) =
w_{\alpha_2}(\alpha_2) = -\alpha_2$. Put $z = u_\gamma \in
U_{\gamma}, u_\gamma \ne 1$. Then $\dot w_c^{-1}z \dot w_c =
u_{-\alpha_2}\in U_{-\alpha_2}$. Further, for every $0< \beta \ne
\alpha_2$ either $\beta + (-\alpha_2)$ is not a root or $\beta
+(-\alpha_2) \in R^+$. Hence $u_\beta u_{-\alpha_2} u_{\beta}^{-1}
= u_{-\alpha_2}v$ for some $v \in U$. Also, for every
$u_{\alpha_2}^\prime \in U_{\alpha_2}$
$$u_{\alpha_2}^\prime
u_{-\alpha_2} u_{\alpha_2}^{\prime-1}\in \langle  U_{-\alpha_2},
U_{\alpha_2}\rangle \,\,\, \text{and} \,\,\,u_{\alpha_2}^\prime
u_{-\alpha_2} u_{\alpha_2}^{\prime-1}\notin B.$$ Recall that $g =
u^{-1}\dot w_c^{-1}$. We may assume $u = v u^\prime_{\alpha_2}$
where the element $v\in U$ does not have factors from
$U_{\alpha_2}$. We have
$$
\begin{array}{ccl}
[g, z] & = & u^{\prime -1}_{\alpha_2}v^{-1}(\dot w_c^{-1}u_\gamma
\dot w_c) v u^\prime_{\alpha_2} u_\gamma^{-1} = u^{\prime
-1}_{\alpha_2}v^{-1}u_{-\alpha_2}v u^\prime_{\alpha_2}
u_\gamma^{-1} \\
& = & u^{\prime -1}_{\alpha_2}u_{-\alpha_2} u^\prime_{\alpha_2}
v^\prime u_\gamma^{-1}
\end{array}
$$
for some $v^\prime \in U$. Put $f = u^{\prime
-1}_{\alpha_2}u_{-\alpha_2}u^{\prime}_{\alpha_2}$ and
$u^{\prime\prime} = v^\prime u_\gamma^{-1}$. We have
$$[g, z] = fu^{\prime\prime}$$ where
$f \in \langle U_{\alpha_2}, U_{-\alpha_2}\rangle, f \notin B$ and
$u^{\prime\prime} \in U.$
\end{proof}

By Lemmas \ref{lem:positive}--\ref{lem:Levi}, we can come up with
the situation when $\Gamma \le G $ corresponds to the root system
generated by $\alpha_1, \alpha_2$ (in our notations), i.e.,
$\Gamma$ is of type $A_2$ (here $\Gamma$ denotes the derived
subgroup of the Levi factor of the corresponding parabolic
subgroup of $G$). By Lemmas \ref{lem:lefttwo} and
\ref{lem:leftonemod}, we have got the following elements in
$\Gamma$ (which are images of commutators of $G$):
$$
v_1 = u_{-\alpha_1}v^\prime,\ v_2 = fv^{\prime\prime},\ u =
u_{\alpha_1}u_{\alpha_2}u^\prime,
$$
where $ 1 \ne u_{-\alpha_1}\in U_{-\alpha_1},\ v^\prime \in
U_\Gamma := \langle U_{\alpha_1}, U_{\alpha_2}\rangle,\ f \in
\langle U_{-\alpha_2}, U_{\alpha_2}\rangle , f \notin B,
\,\,\,v^{\prime\prime}\in U_\Gamma,\,\,\,1\ne  u_{\alpha_1} \in
U_{\alpha_1}, 1\ne \ u_{\alpha_2}\in U_{\alpha_2},\ u^\prime\in
U_{\alpha_1 + \alpha_2}$.

We have to show that the group $\langle v_1,v_2,u\rangle$ is not
solvable.  Consider the groups
$$
P = \langle v_1, u\rangle \,\le\,\tilde{P} = \langle
u_{-\alpha_1}, U_\Gamma \rangle,\,\,\,\,P^\prime = \langle v_2,
u\rangle \,\le\,\tilde{P}^\prime = \langle u_{-\alpha_2},
U_\Gamma\rangle
$$
and the natural homomorphisms $$\theta \colon \tilde{P}
\rightarrow\tilde{P}/R_u,\,\,\,\theta^\prime \colon
\tilde{P}^\prime \rightarrow\tilde{P}^\prime/R_u^\prime$$ where
$R_u$ (resp. $R_u^\prime$) is the unipotent radical of $\tilde{P}$
(resp. $\tilde{P}^\prime$).
 We have
 $$\theta(\tilde P)=\tilde{P}/R_u \cong SL_2(p),
 \,\,\,\theta^\prime(\tilde P^\prime)=
 \tilde{P}^\prime/R_u^\prime \cong SL_2(p) $$
 where $p = 2, 3$. Obviously, $\langle u_{-\alpha_1},
 u_{\alpha_1}\rangle \cong SL_2(p), \langle f,
 u_{\alpha_2}\rangle \cong SL_2(p)$ if $p = 2,3$. Hence
$$\theta(P) = \langle u_{-\alpha_1}, u_{\alpha_1}\rangle \cong SL_2(p),
\,\,\,\theta^\prime(P^\prime) = \langle f, u_{\alpha_2}\rangle
\cong SL_2(p),\,\,\, p= 2, 3.$$

Let us show that $$\Ker\,\theta \cap P \ne
1,\,\,\,\Ker\,\theta^\prime \cap P^\prime \ne 1.$$ Recall that $u$
is regular, so if $\mid K \mid = 2$, then $u^2\in U_{\alpha_1 +
\alpha_2}$, and thus the order of $u$ equals 4. Hence  $u^2 \in
\Ker\theta \cap P\,\, \,\,(u^2 \in \Ker\theta^\prime \cap
P^\prime)$. Let now $\mid K \mid = 3$. Take $h \in P$ (or $h \in
P^\prime$)  such that $\theta(h)$ (or $\theta^\prime(h))$ equals
$\di (-1,-1) \in SL_2(3)$. Then explicit matrix calculations show
that
$$[h,u] = u_{\alpha_2}u_{\alpha_1+\alpha_2} \in \Ker\theta\cap P\,\,\,(\text{or}\,\,\,
[h,u] = u_{\alpha_1}u_{\alpha_1+\alpha_2} \in \Ker\theta^\prime
\cap P^\prime).$$

We proved that $\Ker\,\theta\cap P$ (resp. $\Ker\,\theta^\prime
\cap P^\prime$ )  is not trivial. Let us show that
$\Ker\,\theta\cap P=\Ker\,\theta$ (resp. $\Ker\,\theta^\prime \cap
P^\prime =\Ker\,\theta^\prime $). Note that $\Ker \theta \cong
K^2$ is a 2-dimensional vector $K$-space on which $P$ acts by
conjugation. Since $\theta(P) \cong SL_2(p)$, we have only one
nonzero orbit of $P$ in $\Ker \theta \cong K^2$. Hence $\Ker
\theta \cap P = \Ker\theta \cong K^2$, and therefore $P =
\tilde{P}$. By the same arguments, $\Ker\,\theta^\prime \cap
P^\prime =\Ker\,\theta^\prime $ and $P^\prime = \tilde{P}^\prime$.
Hence
$$P = \langle v_1, u\rangle \,=\,\tilde{P} = \langle u_{-\alpha_1},
U_\Gamma\rangle , \quad P^\prime = \langle v_2, u\rangle
\,=\,\tilde{P}^\prime = \langle u_{-\alpha_2}, U_\Gamma\rangle .$$
Thus, $U_\Gamma, u_{-\alpha_1}, u_{-\alpha_2}$ are all contained
in $\langle v_1, v_2, u\rangle$, and therefore
$$\Gamma = \langle v_1, v_2, u\rangle \cong SL_3(p).$$

{\it Case 4. $\mid X_i\mid = 1$ for every $i$}. Since for all groups
of rank one or two the proposition has been checked, we may assume
$\rank\ G > 2$.

First suppose that the root system corresponding to $G$ does not
contain $D_4$, i.e. is of one of the types $A_r, B_r, C_r, F_4$.
Suppose $X_i = \{\alpha_j\}$, where $j$ is the number of the root
in the standard numbering. Note that by construction of $X$,
neither $\alpha_{j-1}$, nor $\alpha_{j+1}$ belong to $X$. Suppose
that $\alpha_{j+2}\notin X$ or $\alpha_{j-2} \notin X$ (in
particular, this assumption holds if $\alpha_{j\pm 2}$ does not
exist). Then the subgroup $L$ of $G$ generated by
$U_{\pm\alpha_j}$ and $U_{\pm\alpha_{j+1}}$ (or
$U_{\pm\alpha_{j-1}}$) commutes with the elements $U_{\pm \beta}$
for every $\beta \in X \setminus X_i$. Thus we are reduced to the
group $L$ of rank two, and the statement is proved. Let us now
suppose that $\alpha_{j+2} \in X$. Then we can consider the group
$L = \langle U_{\pm \alpha_j}, U_{\pm\alpha_{j+1}},
U_{\pm\alpha_{j+2}}\rangle$ which commutes with the groups
$U_{\pm\beta}$, $\beta \in X\setminus (X_i\cup \{\alpha_{j+2}\})$.
Hence we may assume $\rank\ G = 3$ and $g = \dot w_{\alpha_1}\dot
w_{\alpha_3}u$ for some $u\in U$. Here we have to check the groups
$A_3(p)$, $B_3(p)$, $C_3(p)$, $^2D_4(p)$, $^2A_5(p)$, $^2A_6(p)$,
$p=2, 3$. We can exclude $^2D_4(p)$, $^2A_5(p)$, $^2A_6(p)$, $p=2,
3$, because these groups have a root subgroup $G_\alpha, \alpha =
\alpha_1$ or $\alpha = \alpha_3$, which is isomorphic to
$SL_2(p^2)$, and we can use our considerations for rank one. Since
$A_3(2) \cong A_8$, $B_3(2) \cong C_3(2),$ it remains to calculate
in the groups $A_3(3)$, $B_3(2)$, $B_3(3)$, $C_3(3).$ These groups
are checked by explicit MAGMA calculations. Table 1 contains the
results for the simple groups. For the reductive groups of type
$A_3$ we have $PGL_4(2)\cong PSL_4(2)\cong A_8$, and we only have
to compute $PGL_4(3)$. This group contains no nontrivial 2-radical
elements.

Suppose now that the root system of $G$ is of type $D_r$ or $E_r$.
Let $\beta$ be the root corresponding to the node  with $3$ edges
on the Dynkin diagram. First suppose $\beta \in X$. Then we can
take $\gamma \in \Pi$ which is joined with $\beta$ and disjoint
from all other roots. As $\beta \in X$, we have $\gamma\notin X$,
and $L = \langle U_{\pm\beta}, U_{\pm\gamma}\rangle$ commutes with
every $U_{\pm \delta}, \delta \ne \beta , \delta \in X$. Thus we
may reduce our considerations to groups of rank $2$. Let now
$\beta \notin X$. Suppose $r> 4$. If none of $\alpha_1$,
$\alpha_2$ belongs to $X$, we are reduced to the case of type
$A_2$ treated above. If not, we are reduced to the case of groups
of rank 1. So it remains to consider the case $r=4$, i.e., the
case of the groups $D_4(p), \ p= 2,\ 3.$ This is checked by MAGMA
(see Table 1).
\end{proof}

\section {Groups $^2F_4(q^2)$}

Recall that in light of Remark \ref{rem:st} we have to consider
the groups of type $^2F_4(q^2)$ separately.

If $R$ is a root system and $G_R$ is a connected reductive
algebraic group with root system $R$ defined over some
algebraically closed field, we denote by $\tilde G_R$ the
universal cover of the derived group of $G_R$. If it is clear what
is the root system under consideration, we often drop the
subscript $R$. In particular, throughout this section we denote by
$G$ the twisted Chevalley group $^2F_4(q^2)$, $q =
\sqrt{2^{2m+1}}$, and by $\tilde G$ the simple algebraic group of
type $F_4$ defined over ${\mathbb F}_2$ (identifying it with its
group of ${\overline{\mathbb F}_2}$-points). We have $G\subset
\tilde G$. Correspondingly, tilde always indicates to subgroups of
$\tilde G$. We denote $K={\mathbb F}_{q^2}$.

\begin{theorem} \label{th:F4}
Let $G={}^2F_4(q^2)$. Then $\kappa (G)=2$.
\end{theorem}

\begin{proof}

For $m=0$, the group $G$ is not simple; its derived subgroup (the
Tits group) is checked by MAGMA (see Table 1). So throughout below
we assume $m>0$.

Let $1\ne g \in G$. First suppose $g \in P$ for some parabolic
subgroup $P$. Any parabolic subgroup is conjugate to a standard
parabolic subgroup (see \cite{Ca2}). We may thus assume $P$ to be
a standard parabolic subgroup. We have $P = LV, V = R_u(P)$. We may
assume that the image of $g$ in $P/Z(L)V$ is not trivial (as above)
and reduce the consideration to the group $L/Z(L)$ of semisimple rank
1.

Hence we may assume that $g$ does not belong to any parabolic
subgroup $P$. Then (see \cite[6.4.5]{Ca2}) the order of $C_G(g)$
is prime to $p=2$ (and so is an odd number). Hence $g$ is a
regular semisimple element, and by \cite{Gow} we can get
representatives of any two semisimple conjugacy classes of $G$ in
the form $\sigma =[g, x], \tau = [g, y]$.

Put $H = \langle \sigma, \tau\rangle. $ Suppose $H$ is solvable.
Denote by $I = \{p_1, \dots, p_k\}$ some set of prime divisors of
$\mid H\mid$ and by $H_I$ a Hall subgroup of $H$ corresponding to
$I$. Let $A$ be a maximal normal abelian subgroup of $H_I$.

Let us now consider two separate cases: $m\geq 2$ and $m=1$.

{\it General case $q = \sqrt{2^{2m+1}} ,\,\,\, m \geq 2$.}

\medskip
We have \cite[2.9, p.~76]{Ca2}

$$
\begin{array}{ll}
\mid G \mid = q^{24}(q^2-1)(q^6+1)(q^8-1)(q^{12}+1) =& {}\\
(q^2)^{12}(q^2-1)^2(q^2 +1)^2 ((q^2)^2 +1)^2((q^2)^2- q^2
+1)((q^2)^4-(q^2)^2 + 1), & {}
\end{array}
$$
where $q = \sqrt{2^{2m+1}}$.

\begin{lemma} \label{lem:F4regular}
Let $T$ be a maximal quasisplit torus of $G$. Then there exists $t
\in T$ such that $t$ is a regular element of $\tilde G$, i.e.
$C_{\tilde G}(t)=\tilde S$ is a maximal torus in $\tilde G$.
\end{lemma}

\begin{proof}
Let $\tilde S$ be a maximal torus of $\tilde G$ containing $T$. Let
$\alpha$ be a positive root of $R = F_4$ corresponding to $\tilde
S$, and let $\alpha_T \colon T\rightarrow \overline{\mathbb F}^*_2$
be the restriction of $\alpha$ to $T$.

Let us show that
\begin{equation}
\Imm  \,\alpha_T = K^* \label{eq:41}
\end{equation}
for every $\alpha\in R(F_4)$. We have the following simple root
system
$$\alpha_1 = \epsilon_2- \epsilon_3, \alpha_2 =
\epsilon_3-\epsilon_4, \alpha_3 = \epsilon_4, \alpha_4 =
\frac{1}{2}(\epsilon_1-\epsilon_2-\epsilon_3-\epsilon_4),$$ and
$$T = \langle h_1(t) =
h_{\alpha_1}(t)h_{\alpha_4}(t^\theta),\,\,\,h_2(s)=
h_{\alpha_2}(s)h_{\alpha_3}(s^\theta)\rangle$$ where $s, t \in K^*,
2\theta^2 = 1$. Further,
$$\epsilon_1(h_1(t)) = t^\theta, \epsilon_2(h_1(t)) = t^{1-\theta},
\epsilon_3(h_2(s)) = s, \epsilon_4(h_1(t)) = t^{-\theta}$$ (note
$2(1 -\theta)(1+\theta) = 2 -2\theta^2 = 2-1 =1$),
$$(\epsilon_1 + \epsilon_2)(h_1(t)) = t,
(\epsilon_1 - \epsilon_2)(h_1(t)) = t^{1-2\theta} $$ $\,\,\,\,(
(1-2\theta)(1+ 2\theta) = 1 -4 \theta^2 = 1- 2 =-1)$,
$$(\epsilon_1 \pm \epsilon_3)(h_2(s)) = s^{\pm 1},
(\epsilon_1 \pm \epsilon_4)(h_2(s)) = s^{\pm 1\pm 2\theta},$$
$$(\epsilon_2 \pm \epsilon_3)(h_2(s)) = s^{\pm 1},
(\epsilon_2 \pm \epsilon_4)(h_2(s)) = s^{\pm 1\pm 2\theta}$$
$$(\epsilon_3 + \epsilon_4)(h_2(s)) = s^{2\theta},
(\epsilon_3 - \epsilon_4)(h_2(s)) = s^{2- 2\theta},$$
$$\frac{1}{2}(\epsilon_1\pm\epsilon_2\pm\epsilon_3\pm\epsilon_4)(h_2(s))
 = s^{\pm 1\pm \theta}\,\,\,\text{or}\,\,\,s^{\pm\theta}.$$
Thus we have (\ref{eq:41}). From (\ref{eq:41}) we get
$$\mid \Ker\ \alpha_T \mid \ = (q^2-1)$$
and
\begin{equation}
\mid\bigcup_{\alpha \in R^+(F_4)}\Ker\ \alpha_T \mid \ < (q^2-1)
 \cdot 24 < (q^2-1)^2. \label{eq:42}
\end{equation}

From (\ref{eq:42}) we conclude that the set $ M=T \setminus
\bigcup_{\alpha \in R^+(F_4)}\Ker \ \alpha_T $ is not empty. Any
element $t\in M$ is regular. The lemma is proved.
\end{proof}

\begin{lemma} \label{lem:F4-p}
There exists a prime $p \ne 2, 3, (p,q^2-1) = 1$ such that $p\mid
q^2+1$ or $p \mid q^4 +1$.
\end{lemma}

\begin{proof}
This follows from the fact that $(q^2 - 1, q^2+ 1) = 1$ and $(q^2
+1, q^4 +1) = 1$.
\end{proof}

\begin{lemma} \label{lem:maxtorus}
Let $R$ be a root system, and let $G_R$ be a connected reductive
group. Further, let $A\subseteq G_R$ be a finite abelian subgroup
consisting of semisimple elements and such that $(\mid A\mid ,
\mid W(R)\mid) = 1$. Then there exists a maximal torus $S$ in
$G_R$ such that $A \subseteq S$.
\end{lemma}

\begin{proof}
Let  $G_R = S^\prime G^\prime_{R}$, where $S^\prime\le Z(G_R)$ is
a torus of $G_R$ and $G_{R}^\prime$ is semisimple. Hence $Z(G_R) =
S^\prime A^\prime$, where $A^\prime = Z(G_R^\prime)$ is a finite
abelian group. Suppose $A \subseteq Z(G_R)$. Since $(\mid A\mid ,
\mid W(R)\mid) = 1$, we have $(\mid A\mid , \mid A^\prime\mid) =
1$ (because $\mid W(R)\mid$ is divisible by $\mid A^\prime\mid$),
and hence $A \le S^\prime$. Suppose $a \notin Z(G)$ for some $a
\in A$. Let $S$ be a maximal torus of $G_R$ containing $a$. By
\cite[Theorem 3.5.3]{Ca2}, we have
$$C_{G_R}(a) = \langle S, U_{ \alpha},
\dot w\,\,\,\mid\,\,\,\alpha(a) = 1, w \in C_{W(R)}(a)\rangle,$$
$$C_G(a)^0 = \langle T_1, U_{ \alpha},
\,\,\,\mid\,\,\,\alpha(a) = 1\rangle.$$ Hence $\mid
C_G(a)/C_G(a)^0\mid$ divides $\mid W(R)\mid$, and therefore $ A\le
C_G(a)^0 \ne G_R$. To finish the proof, we use induction by $\mid
R\mid$.
\end{proof}

\bigskip

Before going over to the proof of the assertion of the theorem, we
shall describe some general construction (parallel to that of Lemma
\ref{lem:normalizer}).

Let $G_R$ be a connected semisimple group corresponding to a root
system $R$, and let $S$ be a maximal torus of $G_R$. Further, let
$M \subseteq S$, let $g \in N_{G_R}(M)$, and let $g = u \dot w v $
be a Bruhat decomposition of $g$ in $G_R$ with respect to a Borel
subgroup containing $S$. We may assume $\dot w v \dot w^{-1} \in
U^-$. Let $s \in M$. Then
$$g s g^{-1} = u \dot w v s v^{-1} \dot w^{-1} u^{-1} = u w(s)
v^\prime u^{-1} = s^\prime \in M \subseteq S,$$ where $v^\prime\in
U^-$. Hence $w(s)v^\prime = u^{-1}s^\prime u = s^\prime[s^{\prime
-1}, u^{-1}]$. Since $[s^{\prime -1}, u^{-1}] \in U, v^\prime = [
w(s)^{-1},v] \in U^-$, we have $[s^{\prime -1}, u^{-1}] = 1,[
w(s)^{-1},v] = 1, s^\prime = w(s)$. Since we can consider any $s
\in M$, we have $u, v \in C_{G_R}(M)$. Now we have a homomorphism
$$\phi\colon N_{G(R)}(M)\rightarrow W(R)$$
with
\begin{equation}
\Ker\phi = C_{G(R)}(M)^0.\label{eq:43}
\end{equation}

\bigskip

We can now go over to the proof of Theorem \ref{th:F4}.

Set $\sigma = t$, where $t$ is chosen as in Lemma
\ref{lem:regular}. Let $\tau$ be an element of order $p$ (it
exists by Lemma \ref{lem:F4-p}).  Denote by $I$ the set consisting
of $p$ and all prime divisors of $q-1$. Since all Hall subgroups
$H_I$ are conjugate and each element of order $p$ belongs to one
of those, we may assume $t \in H_I$ and some element $\tau^\prime$
of order $p$ is also in $H_I$.

Note that $2,3 \nmid (q^2-1)$. Since $\mid W(F_4)\mid = 2^73^2$,
we have $\sigma, \tau^\prime \in C_{\tilde G}(A)^0$ (by
(\ref{eq:43}) and Lemma \ref{lem:maxtorus}). Then $A \subseteq T
\subseteq \tilde{T}$ where $\tilde{T}$ is the unique maximal torus
of $\tilde G$ containing $T$ (recall that $T$ contains a regular
semisimple element of $\tilde G$).

Denote by $R\subset R(F_4)$ the minimal (with respect to
inclusion) root subsystem such that
$$\sigma, \tau^\prime \in G_R = \langle \tilde{T}, U_\alpha\,\,\,\mid
\,\,\,\alpha\in R\rangle.$$ First note that $R\neq R(F_4)$ because
otherwise we would have $A\subseteq Z(F_4)=1$ (recall that $H \leq
C_{\tilde{G}}(A)^0$). Second, note that $R \ne \emptyset $ because
$\tau^\prime \notin T =\tilde{T}^F$. Set $G_R^\prime = \langle
U_\alpha\,\,\,\mid \,\,\,\alpha\in R\rangle$. Then $ G_R =
SG_R^\prime$ where $S\le \tilde{T}\cap Z( G_R)$ is a subtorus of
$\tilde{T}$. Then $Z( G_R) = SZ(G_R^\prime)$. Since the orders of
$\sigma, \tau^\prime$ are prime to $2,3$, we have $\sigma,
\tau^\prime \notin Z( G_R^\prime)$, and hence so are the orders of
their images $\bar{\sigma}, \bar{\tau}^\prime$ in
$\overline{G}_R^{\prime} = G_R^\prime/(Z(G_R^\prime)\cap S)$. Now
we have a semisimple group $\overline{G}^\prime_R$ with a maximal
torus $\overline{T} = \tilde{T}/S$ which contains the solvable
group $\overline{H}_I = \langle \bar{\sigma},
\bar{\tau}^\prime\rangle \ne 1 $, where $\bar{\sigma}\in
\overline{T}$ is a regular element. Let $A_1$ be a maximal abelian
normal subgroup of $\overline{H}_I$. Then $A_1 \subseteq
\overline{T}$ and $A_1 \nsubseteq Z(\overline{G}^\prime_R)$ (note
that $2, 3$ are the only primes dividing both $\mid W(R)\mid$ and
$Z(\overline{G}_R^\prime)$). By (\ref{eq:43}), we have
$$\bar{\sigma}, \bar{\tau}^\prime \in C_{\overline{G}_R^\prime}(A_1)^0 = \langle
\overline{T}^\prime, U_\beta \,\,\mid \,\,\beta(A_1) =1\rangle =
\langle \overline{T}^\prime, U_\beta \,\,\mid \,\,\beta \in
R^\prime \subsetneqq R\rangle.$$ Hence
$$\sigma, \tau^\prime \in \langle \tilde{T},
U_\beta \,\,\mid \,\,\beta \in R^\prime\rangle.$$ This is a
contradiction with the choice of $R$.

\bigskip

Let us now consider the last remaining special case.

{\it Case $q = \sqrt{2^3}$}.

\bigskip
Here $\mid G\mid = 2^{36}\,\cdot 3^5\, \cdot 5^2\, \cdot
7^2\,\cdot 13^2\,\cdot 19\,\cdot 37\cdot 109 $.

\bigskip

Let $\mid \langle \sigma\rangle \mid = 109,\,\,\, \mid\langle \tau
\rangle \mid = 37$, and let $H_0 \subseteq H$ be a Hall subgroup of
$H$ of order $37\cdot 109$. Since $(37, 109 -1) = 1$, the group $H_0
= \langle h\rangle$ is cyclic of order $37\cdot 109$.

Let, as above,  $\tilde{G}$ denote the simple algebraic group of
type $F_4$ over the field ${\mathbb F}_2$, and let $F$ be the
Frobenius map of $G$ such that $G = \tilde{G}^F$. Since $h \in
\tilde{G}^F$, the centralizer $C_{\tilde{G}}(h)$ is an $F$-stable
connected reductive group (\cite[3.5.6]{Ca2}) which, in turn,
contains an $F$-stable maximal torus $\tilde T$ (which is also a
maximal torus of $\tilde{G}$). Hence $h \in \tilde T^F$. But
$$\mid \tilde T^F\mid = \prod_{i=1}^4( q- \epsilon_i)$$ where each
$\epsilon_i$ is a root of unity \cite[3.3.5]{Ca2}. Since
$$\mid q - \epsilon_i\mid \leq q + 1 = \sqrt{8} +1 \leq 4 ,$$
we conclude that $\mid \tilde T^F\mid \leq 256 < 37\cdot 109 $.
Contradiction.

\medskip

The theorem is proved.
\end{proof}

\section{Groups generated by 3-transpositions} \label{sec:transp}

In this section we show that the estimate of Proposition \ref{prop:bigrank} is sharp
as follows from the case of
groups generated  by 3-transpositions (see \cite{Fi}, \cite{As} for
definitions and notations).

\begin{defn}\cite{Fi}
Let $G$ be a finite group generated by a class $D$ of conjugate
involutions such that any pair of non-commuting elements of $D$
generates a dihedral group of order $6$; then $D$ is a class of
{\emph{conjugate $3$-transpositions}} of $G$.
\end{defn}

Equivalently, the product of any two involutions from $D$ is of order
1, 2, or 3.

\begin{prop}
Let $G$ be a finite group generated by a class $D$ of conjugate
$3$-transpositions. Then any element of $D$ is $2$-radical.
\end{prop}

\begin{proof} This is an immediate consequence of
\cite[Cor.~1.6]{Fi}.
\end{proof}

\begin{corol}
Let $G$ be one of the following groups:
\begin{itemize}
\item a symmetric group $S_n$;
\item a symplectic group $\Sp (2n,2) (n\geq 2)$;
\item an orthogonal group $\Or^{\mu}(2n,2)$ for $\mu\in\{-1,1\}$ and $n\geq 2$;
\item a unitary group $\PSU (n,2) (n\geq 4)$;
\item an orthogonal group $\Or^{\mu ,\pi}(n,3)$ for $\mu\in\{-1,1\}$,
$\pi\in\{-1,1\}$, and $n\geq 4$;
\item one of Fischer's groups $Fi_{22}$, $Fi_{23}$, $Fi_{24}$.
\end{itemize}
Then $G$ contains a nontrivial $2$-radical element.
\end{corol}

\begin{proof} This immediately follows from the above proposition taking into
account the fact that all the listed groups are generated by a class
of conjugate 3-transpositions \cite{Fi}.
\end{proof}

\section{ Sporadic groups}
\label{sec:spor}


\begin{prop} \label{prop:spor}
Let $G$ be a sporadic simple group. Then $\kappa (G) = 3$ for
$G=Fi_{22}, Fi_{23}$ and $\kappa (G)=2$ for all the remaining groups.
\end{prop}

More precisely, we shall prove that if $g\ne 1$ is a 2-radical
element of a sporadic simple group $G$, then $G=Fi_{22}$ or
$G=Fi_{23}$ and $g$ is a 3-transposition. (In the latter cases MAGMA
computations show that $g$ is not a 3-radical element.)

The proof goes case by case. Apart from the theoretical arguments
presented below, we used MAGMA for rechecking them (in all the cases
except for the Monster). For larger sporadic groups we had to
replace most standard MAGMA procedures with our own ones in order to
avoid storing the whole group and large subgroups. In particular, to
check whether a subgroup under consideration is not solvable, we
used the Hall--Thompson criterion \cite{Th}: a group $H$ is
nonsolvable if and only if it contains nonidentity elements $a, b,
c$ of pairwise coprime orders such that $abc=1$.

Both in the theoretical proof and in the computer-aided one, we rely
on the ATLAS classification of conjugacy classes of maximal cyclic
subgroups \cite{Atlas}.

Let us now prove the proposition. The exposition below is sometimes
sketchy, we omit some cases where the proof uses arguments similar
to earlier ones.

The main idea is very simple. We first consider the elements of
prime orders. It turns out that in most cases one can include a
given element $g$ of prime order $p$ of a group $G$ in its proper
simple subgroup $H$. If there is a single conjugacy class of cyclic
subgroups of order $p$, it is enough to indicate $H$ whose order is
divisible by $p$. In the case where there are several conjugacy
classes of cyclic subgroups of order $p$, more subtle arguments are
needed. We either use ATLAS information on elements $h$ of order
$mp$ for some $m$ whose powering gives $g$ and try to include $h$ in
some proper simple subgroup $H$, or use some information on subgroup
structure of $G$ from the literature. Finally, if $g$ is not
contained in any proper simple subgroup of $G$, it happens that its
normalizer $N=N_G(g)$ is the unique maximal subgroup of $G$
containing $g$. In that case, one can take $x\in N$ and get
$a=[g,x]\in \left<g\right>$, and take $y$ such that $b=[g,y]\notin
N$. Then $\left<a,b\right>=G$ is not solvable.

If an element $g$ under consideration is of composite order $mp$, we
note that it belongs to the centralizer of $h=g^m$ which is of prime
order $p$. It remains to use the information from ATLAS on the
centralizers of elements of prime orders in sporadic groups. It
turns out that in many cases the structure of $C_G(h)$ is as
follows: it contains a normal subgroup $Z$ of small exponent such
that the quotient $G'=C_G(H)/Z$ is either a smaller simple group or
an extension of a simple group by a group of small exponent. Thus if
$g$ is a $2$-radical element of sufficiently large exponent, then
its image in $G'$ is a nonidentity 2-radical element, and we arrive
at a contradiction by induction. In some cases, elements of small
exponents require separate consideration.

Below we mostly present theoretical arguments as above for elements
of prime orders. We present a more detailed proof for the
baby-monster $B$ and a complete proof for the monster $M$.

We follow the subdivision of sporadic groups from ATLAS.

\medskip
{\it Mathieu groups}: $M_{11}, M_{12}, M_{22}, M_{23}, M_{24}.$

\medskip
$M_{11}$. The elements of orders 2, 3 and 5 are included in $A_5$,
and of order 11 --- in $PSL_2(11)$.

$M_{12}$. Any element of order 11 is included in $PSL_2(11)$. All
the remaining ones, of types 2A, 2B, 3A, 3B, 5A, are included in
$A_5$ (according to \cite[p.~33]{Atlas-book}, $M_{12}$ contains
$A_5$'s of types (2A, 3B, 5A) and (2B, 3A, 5A)).

$M_{22}$. The elements of orders 2, 3, 5 and 7 are included in
$A_7$, and of order 11 --- in $PSL_2(11)$.

$M_{23}$. The elements of orders 2, 3, 5 and 11 are included in
$M_{11}$, of order 7 --- in $A_7$, and the normalizer $N=23\cdot 11$
of an element $g$ of order 23 is the unique maximal subgroup of
$M_{23}$  containing  $g$, so we can apply the argument mentioned
above.

$M_{24}$. Any element of order 23 is included in $PSL_2(23)$, of
order 11 --- in $M_{23}$, of order 7 --- in $PSL_2(7)$, and of order
5 --- in $A_5$. According to \cite[p.~96]{Atlas-book}, $M_{24}$
contains $A_5$'s of types (2B, 3A, 5A) and (2B, 3B, 5A), so it
remains to consider the class 2A. Fix an element $z$ of type 2B. We
have $C_G(z)=E_{2^6}\cdot S_5$, where $E_{2^6}$ is an elementary
abelian subgroup. Let $g$ be any involution of $A_5\subset S_5$.
Since $g$ centralizes $z$, it cannot be conjugate to $z$, hence $g$
is of type 2A, and we are done.

\medskip
{\it Leech lattice groups}: $HS, J_2, Co_1, Co_2, Co_3, McL, Suz$.
\medskip

Here we shall be a little sketchy describing only the largest Conway
group $Co_1$ among the three ones.

$HS$. The elements of orders 3, 7 and 11 are included in $M_{22}$.
According to \cite[p.~80]{Atlas-book}, there is an $M_{11}\subset
HS$ containing elements of types 2A and 5C, and there is an $A_5$
containing elements of type 2B and 5A. The remaining class 5B also
has a representative lying in $A_5$ \cite[p.~274]{GLS}.

$J_2$. Any element of order 7 can be included in $PSL_3(2)$.
According to \cite[p.~42]{Atlas-book}, there are $A_5$'s of types
(2B, 3A, 5CD), (2A, 3B, 5AB), thus including the elements of all the
other classes.

$McL$. There are no problems with the elements of orders 2, 7 and 11
--- they can all be included, say, in $M_{11}$. By
\cite[p.~100]{Atlas-book}, there is a subgroup $PSU_3(5^2)$
containing representatives of 3B, 5A and 5B. It remains to consider
the class 3A. Take an element of order 9 in $PSU_4(3^2)$. According
to \cite{Atlas}, its cube belongs to 3A.

$Suz$. Any element of order 13 belongs to a maximal subgroup
$G_2(4)$, and hence to an even smaller subgroup $PSL_2(13)$. The
elements of orders 7 and 11 belong to $M_{11}$. On
\cite[p.~131]{Atlas-book} we find an $A_7$ containing
representatives of 2B, 3C and 5B, a $PSL_3(3)$ containing
representatives of  3B, and a $PSL_2(25)$ containing representatives
of 5A and 5B. It thus remains to consider the classes 2A and 3A. To
treat 2A, take an element of order 8 in $M_{11}$, then its cube is
of type 2A \cite{Atlas}. Similarly, the fifth power of an element of
order 15 in $J_2$ is of type 3A.

Conway groups: we shall skip the arguments for $Co_2$, $Co_3$.

$Co_1$. The elements of orders 23 and 11 belong to $M_{23}$, and
those of order 13 --- to $Suz$. The classification of $A_5$'s
\cite{Wi83} gives subgroups of types (2B, 3A, 5A), (2C, 3A, 5B),
(2C, 3B, 5C), (2B, 3B, 5A), (2B, 3A, 5A).   According to \cite{Cu},
the classes 7A and 7B have their representatives in $A_7$ and
$PSL_2(7)$, and the class 3D, as 3A, belongs to $A_5$. It remains to
consider 2A. One can take an element of order 18 in $Co_3$, its 9th
power is of type 2A.

\medskip
{\it Monster sections}: $He, HN, Th, Fi_{22}, Fi_{23}, Fi'_{24}, B,
M$.
\medskip

Here we shall skip $HN$ (which can be treated using
\cite[p.~166]{Atlas-book} and \cite{NW}) and two larger Fischer
groups.

$He$. The elements of order 17 belong to $PSp_4(4)$, and hence to
$PSL_2(16)$. The elements of order 5 lie in $A_5$. We have to
consider the classes 2A, 2B, 3A, 3B, 7A, 7C and 7D (7B is a power of
7A and 7E is a power of 7D). First we use the information on
(2,3,7)-subgroups from \cite[p.~104]{Atlas-book}: a subgroup of type
(2A, 3B, 7C) is contained in $7:3\times PSL_3(2)$ (and hence 2A
belongs to $PSL_3(2)$), and a subgroup of type (2B, 3A, 7AB) is
contained in $S_4\times PSL_3(2)$ (and hence 7A belongs to
$PSL_3(2)$ too). Next, we use the information on the centralizers of
involutions \cite[p.~277]{GLS}. Since 7D and 3B commute with 2B,
they both belong to $PSL_3(2)$. Since 3A commutes with 2A, it
belongs to the centralizer of 2A, and hence to $PSL_3(4)$. As to 2A
and 2B, the same argument as in the case $M_{24}$ applies, and we
conclude that 2A belongs to $PSL_3(2)$ and 2B belongs to $PSL_3(4)$.
Finally, since 7C commutes with 3A, it belongs to the centralizer of
3A and hence to $A_7$.

$Th$. The normalizer $N=31\cdot 15$ of an element $g$ of order 31 is
the unique maximal subgroup of $Th$  containing  $g$, so we can
proceed as in the case of an element of order 23 in $M_{23}$. Any
element of order 19 belongs to $PSL_2(19)$, of order 13 --- to
${}^3D_4(2)$, of orders 2, 5 and 7 --- to $A_7$. It remains to treat
three classes of elements of order 3. Take an element of order 21 in
$PSL_5(2)$, its 7th power is of type 3A. Taking elements of orders 9
and 15 in $2^{1+8}\cdot A_9$, we obtain 3B and 3C as their 3rd and
5th power, respectively.

$Fi_{22}$. First recall that this group does contain 2-radical
elements, namely, those of the class 2A (3-transpositions), see
Section \ref{sec:transp} above. Any element of order 13 belongs to
$O_7(3)$, and hence to  $PSL_3(3)$. The elements of orders 5, 7, 11
lie in $M_{22}$. We have to consider the classes 2B, 2C, 3A, 3B, 3C,
3D. According to \cite[p.~163]{Atlas-book},  there is an $M_{12}$
containing representatives of 2B, 2C, 3C, 3D. We include 3A in
$A_{10}$ representing it as the 5th power of an element of order 15
in $A_{10}$. Similarly, we represent 3B as the 6th power of an
element of order 18 in $O^+_8(2)$.

$B$. The normalizer $N=47\cdot 23$ of an element $g$ of order 47 is
the unique maximal subgroup of $B$  containing  $g$, so we can
proceed as above. The cases of elements of orders 31, 23, 19, 17,
13, 11 and 7 are easy: those of order 31 belong to $PSL_2(31)$, of
order 19 --- to $Th$, and all the remaining ones can be included,
say, in $Fi_{23}$. Furthermore, we use the classification of $A_5$'s
\cite[Theorems 5.1, 5.2]{Wi93}: in particular, there are subgroups
of types (2B, 3A, 5A), (2D, 3B, 5B) and also those containing 2C. It
remains to consider 2A. We get it as the 13th power of an element of
order 26 in $Fi_{23}$.

Let now $g$ be an element of composite order $mp$, $m\ge p$. As
$p\le 5$, it suffices to use information on the centralizers of
the elements of orders 2, 3 and 5. We have $C_B($2A$)=2\cdot
({}^2E_6(2)):2$, $C_B($2B$)=2_+^{1+22}\cdot Co_2$,
$C_B($2C$)=(2^2\cdot F_4(2)):2$, $C_B($2D$)=2^9\cdot 2^{16}\cdot
O_8^+(2)\cdot 2$, $C_B($3A$)=3\times Fi_{22}:2$,
$C_B($3B$)=3_+^{1+8}:2_-^{1+6}\cdot PSU_4(2)$, $C_B($5A$)=5\times
HS:2$, $C_B($5B$)=5_+^{1+4}:2_-^{1+4}\cdot A_5$.

First suppose $g$ is of odd order $mp$, $m > p$. If $p=3$, then
$g$ centralizes either 3A or 3B. As the exponent of the
extraspecial group $3_+^{1+8}$ equals 3, we get the image of $g$
of order at least 5 in either $Fi_{22}$ or $PSU_4(2)$ whose
2-radical elements can only be of order 2 or 3. Thus $g$ is not
2-radical. (Note that this argument does not work for the elements
of order 9 which will be considered separately.) If $p=5$, we have
to consider the elements of orders 35 and 55 which all centralize
5A. Hence each of them maps to a nonidentity element of $HS$, and
we are done. The elements of order 25 centralize 5B. As the
exponent of the extraspecial group $5_+^{1+4}$ equals 5, each of
them maps to a nonidentity element of $A_5$ which cannot be
2-radical. To finish with the case of odd order, it remains to
consider the elements of order 9. According to \cite{Atlas}, both
9A and 9B can be represented as the 4th power of an element of
order 36. Hence any element of order 9 centralizes either 2B or 2D
and thus belongs to either $Co_2$ or $O_8^+(2)$, and we are done.

Suppose now $g$ is of even order $2m$ so that $g$ centralizes an
involution of $B$. If $m$ is odd, then the image of $g$ in the
simple group involved in the centralizer of the corresponding
involution is nonidentity, and we are done. So assume $m$ to be
even, i.e. $g$ is of order $4n$. The elements of order 4 were
checked by MAGMA, so suppose $n>1$. According to \cite{Atlas},
there are no elements of order $4n$, $n>1$, powering to 2A. If $g$
centralizes 2B, then it maps to a nonidentity element of $Co_2$,
and we are done. According to \cite{Atlas}, the elements of order
$4n$, $n>1$, powering to 2C are 12T, 20H and 52A, they were
checked separately by MAGMA. Finally, suppose that $g$ centralizes
2D. If $n>2$, then taking into account that $C_B(2D)<2^9\cdot
2^{16}\cdot PSp_8(2)$, we conclude that $g$ maps to a nonidentity
2-radical element of order greater than 2 in $PSp_8(2)$ which
contradicts to MAGMA computations in that group. Thus it remains
to check the elements of order 8 powering to 2D, i.e. 8G, 8J, 8K,
8M and 8N. This was also done by MAGMA.

$M$. In this case no additional MAGMA computations were needed, we
only used the results for smaller groups. Our approach mimics the
case of the baby-monster.

The normalizer $N=41\cdot 40$ of an element $g$ of order 41 is the
unique maximal subgroup of $M$ containing $g$, so we can proceed
as above. Relying on the existing information on maximal subgroups
of $M$ \cite{BrW}, we include the elements of orders 71, 59, 47,
31, 29, 23, 19, 17, 11 in $PSL_2(71)$, $PSL_2(59)$, $B$, $B$,
$Fi'_{24}$, $B$, $B$, $B$, $B$, respectively. Representatives of
all the remaining classes appear in \cite{No}: Table 1 on p.~201
gives 13A and 13B lying in $PSL_3(3)$, in Section 5 there are
exhibited 7A and 7B lying in $PSL_3(2)$, and the list of $A_5$'s
in Table 3 on p.~202 contains representatives of all classes of
elements of orders 2, 3 and 5.

Let now $g$ be an element of composite order $mp$, $m\ge p$. Our
arguments are similar to the previous case. As for $B$, we have
$p\le 5$, and it suffices to use information on the centralizers
of the elements of orders 2, 3 and 5. We have $C_M($2A$)=2\cdot
B$, $C_M($2B$)=2_+^{1+24}\cdot Co_1$, $C_M($3A$)=3\times
Fi'_{24}$, $C_M($3B$)=3_+^{1+12}.2Suz$, $C_M($3C$)=3\times Th$,
$C_M($5A$)=5\times HN$, $C_M($5B$)=5_+^{1+6}:2J_2$.

First suppose $g$ is of odd order $mp$, $m \ge p$. If $p=3$, then
$g$ centralizes either 3A, or 3B, or 3C. As the exponent of the
extraspecial group $3_+^{1+12}$ equals 3, we get the image of $g$
of order at least 5 in either $Fi'_{22}$, or $Suz$, or $Th$ which
do not contain 2-radical elements. Thus $g$ is not 2-radical. If
$p=5$, we have to consider the elements of orders 25, 35, 45, 55,
95 and 105. Any of those centralizes either 5A or 5B and hence
maps to a nonidentity element of either $HN$ or $J_2$. (We use the
fact that the exponent of the extraspecial group $5_+^{1+6}$
equals 5.)

If $g$ is of even order $2m$, it centralizes either 2A or 2B. If
$m>2$, then $g$ maps to a nonidentity element of either $B$ or
$Co_1$ which do not contain 2-radical elements. Thus $g$ is not
2-radical and we are done. Let now $g$ be of order 4. Any
4A-element is the 11th power of 44A and hence belongs to $B$. The
square of a 4B-element belongs to 2A \cite{Atlas}. Therefore 4B
centralizes 2A and thus maps to a nonidentity element of $B$.
According to \cite{Atlas}, the 4th power of any element of order
16 belongs to 4C, hence 4C lies, say, in $Fi'_{24}$. Finally, 4D
is the cube of 12J whose 4th power is 3C. Therefore 12J
centralizes 3C and hence so does 4D. Thus 4D belongs to $Th$, and
we are done.

\medskip
{\it Pariahs}: $J_1, J_3, J_4, Ru, O'N, Ly$.
\medskip

$J_1$. The normalizer $N=19\cdot 6$ of an element $g$ of order 19 is
the unique maximal subgroup of $J_1$  containing  $g$, and the above
argument applies. If the order of $g$ equals 7, its normalizer $N$
equals $7\cdot 6$ and is also a maximal subgroup of $J_1$ but is
contained in another maximal subgroup of order 168. However, taking
$x\in N$ and $y$ of order 3, we get $a=[g,x] \in \left<g\right >$
and $b=[g,y]$ of order 15. Since $b$ is outside of both above
mentioned maximal subgroups, we have $\left<a,b\right>=J_1$. The
elements of order 11 belong to $PSL_2(11)$, and the elements of
orders 2, 3 and 5 belong to $A_5$.

$J_3$. The elements of orders 2 and 5 belong to $A_5$, those of
orders 17 and 19 belong to $PSL_2(17)$ and $PSL_2(19)$,
respectively. Taking an element of order 9 in $PSL_2(17)$, we obtain
3B as its cube, and taking an element of order 15 in $PSL_2(16)$, we
obtain 3A as its 5th power.

$J_4$. For $p=43$ or 29, the normalizer of $g$ of order $p$ is the
unique maximal subgroup containing $g$, and we apply the above
argument. The elements of order 37 lie in $PSU_3(11^2)$, of order 31
--- in $PSL_2(32)$, and of orders 3, 5, 7, and 23 --- in $M_{24}$.
It remains to consider the classes 2A, 2B, 11A, 11B. The
centralizers of each of 2A and 2B contain $M_{22}$, and we embed
both 2A and 2B in $M_{22}$ using the same argument as in the case
$M_{24}$ above. According to \cite[Prop.~22 and Prop.~26]{J}, we
have 11A$\in C$(2B) and 11B$\in C$(2A), so they are both included in
$M_{22}$ too.

$Ru$. The elements of orders 29 and 13 lie in the corresponding
$PSL$'s, and those of orders 7 and 3 lie in $A_7$. The information
on alternating subgroups in \cite[p.~126]{Atlas-book} gives 2B, 5A
and 5B contained there. 2A appears as the square of an element of
order 4 in $A_6$.

$O'N$. The elements of order 31 lie in $PSL_2(31)$, of order 19
--- in $PSL_3(7)$, and of orders 11, 5, 3 and 2 --- in $M_{11}$.
As to the classes 7A and 7B, the first appears as the square of an
element of order 14 in $PSL_3(7)$, and the second belongs to
$PSL_2(7)$ \cite[Section 4, p.~471]{Wi85}.

$Ly$. For $p=67$ and 37 we use the same maximal subgroup argument as
above. The elements of order 31 belong to $G_2(5)$, and hence to
$PSL_3(5)$, and those of orders 11, 7 and 2 --- to $A_{11}$. The
classification of $A_5$'s \cite[Section 6, p.~407]{Wi84} shows that
3B and 5B are included in $A_5$. Both 3A and 5A lie in $G_2(5)$:
they can be obtained as the 3th power of an element of order 9 and
the 4th power of an element of order 20, respectively.

\medskip

To finish the proof of the proposition, it remains to check all
small groups of Lie type appearing in the above arguments. This was
done by straightforward computations.

Proposition \ref{prop:spor}, and hence Theorems \ref{th:lie1} and
\ref{main:solv}, are proved. \qed



\end{document}